\documentclass[12pt,sort&compress]{elsarticle}

\usepackage[pdfborder={0 0 0},colorlinks,allcolors=blue]{hyperref}

\usepackage{mystyle}
\usepackage{afterpage}

\theoremstyle{definition}
\newtheorem{remark}{Remark}

\journal{}

\begin{document}

\begin{frontmatter}

\title{Interpolation-based immersed finite element and isogeometric analysis}

\author[ucsd]{Jennifer E. Fromm}
\author[colorado]{Nils Wunsch}
\author[ucsd]{Ru Xiang}
\author[ucsd]{Han Zhao}
\author[colorado]{Kurt Maute}
\author[colorado]{John A. Evans} \author[ucsd]{David Kamensky\corref{cor1}}
\ead{dmkamensky@eng.ucsd.edu}
\cortext[cor1]{Corresponding author}

\address[ucsd]{University of California San Diego, San Diego, CA 92093, USA}
\address[colorado]{University of Colorado Boulder, Boulder, CO 80309, USA}

\begin{abstract}
We introduce a new paradigm for immersed finite element and isogeometric methods based on interpolating function spaces from an unfitted background mesh into Lagrange finite element spaces defined on a foreground mesh that captures the domain geometry  but is otherwise subject to minimal constraints on element quality or connectivity.  This is a generalization of the concept of Lagrange extraction from the isogeometric analysis literature and also related to certain variants of the finite cell and material point methods.  Crucially, the interpolation may be approximate without sacrificing high-order convergence rates, which distinguishes the present method from existing finite cell, CutFEM, and immersogeometric approaches.  The interpolation paradigm also permits non-invasive reuse of existing finite element software for immersed analysis.  We analyze the properties of the interpolation-based immersed paradigm for a model problem and implement it on top of the open-source FEniCS finite element software, to apply it to a variety of problems in fluid, solid, and structural mechanics where we demonstrate high-order accuracy and applicability to practical geometries like trimmed spline patches.  
\end{abstract}

\begin{keyword}
immersed finite element method \sep
immersogeometric analysis \sep 
trimmed isogeometric analysis  \sep
CutIGA  \sep
Lagrange extraction
\end{keyword}

\end{frontmatter}

\section{Introduction}
\label{sec:introduction}

\begin{figure}[t]
	\begin{center}
	\def\svgwidth{13.0cm}
\begingroup%
  \makeatletter%
  \providecommand\color[2][]{%
    \errmessage{(Inkscape) Color is used for the text in Inkscape, but the package 'color.sty' is not loaded}%
    \renewcommand\color[2][]{}%
  }%
  \providecommand\transparent[1]{%
    \errmessage{(Inkscape) Transparency is used (non-zero) for the text in Inkscape, but the package 'transparent.sty' is not loaded}%
    \renewcommand\transparent[1]{}%
  }%
  \providecommand\rotatebox[2]{#2}%
  \newcommand*\fsize{\dimexpr\f@size pt\relax}%
  \newcommand*\lineheight[1]{\fontsize{\fsize}{#1\fsize}\selectfont}%
  \ifx\svgwidth\undefined%
    \setlength{\unitlength}{262.28303444bp}%
    \ifx\svgscale\undefined%
      \relax%
    \else%
      \setlength{\unitlength}{\unitlength * \real{\svgscale}}%
    \fi%
  \else%
    \setlength{\unitlength}{\svgwidth}%
  \fi%
  \global\let\svgwidth\undefined%
  \global\let\svgscale\undefined%
  \makeatother%
  \begin{picture}(1,0.32782642)%
    \lineheight{1}%
    \setlength\tabcolsep{0pt}%
    \put(0,0){\includegraphics[width=\unitlength,page=1]{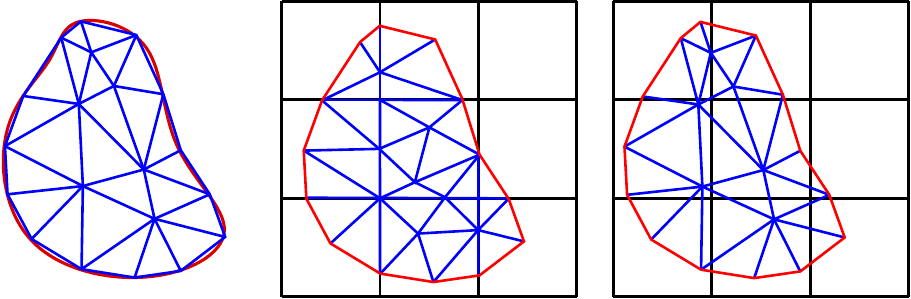}}%
    \put(0.19280988,0.27729961){\color[rgb]{0,0,0}\makebox(0,0)[lt]{\lineheight{1.25}\smash{\begin{tabular}[t]{l}(a)\end{tabular}}}}%
    \put(0.55840528,0.27102039){\color[rgb]{0,0,0}\makebox(0,0)[lt]{\lineheight{1.25}\smash{\begin{tabular}[t]{l}(b)\end{tabular}}}}%
    \put(0.92597178,0.27102039){\color[rgb]{0,0,0}\makebox(0,0)[lt]{\lineheight{1.25}\smash{\begin{tabular}[t]{l}(c)\end{tabular}}}}%
  \end{picture}%
\endgroup%

	\caption{(a) Domain discretized using a body-fitted mesh; (b) Domain discretized using a tensor-product background mesh with a background-fitted foreground mesh; (c) Domain discretized using a tensor-grid background mesh with a background-unfitted foreground mesh.} 
	\label{fig:immersed_meshes}
    \end{center}
\end{figure}
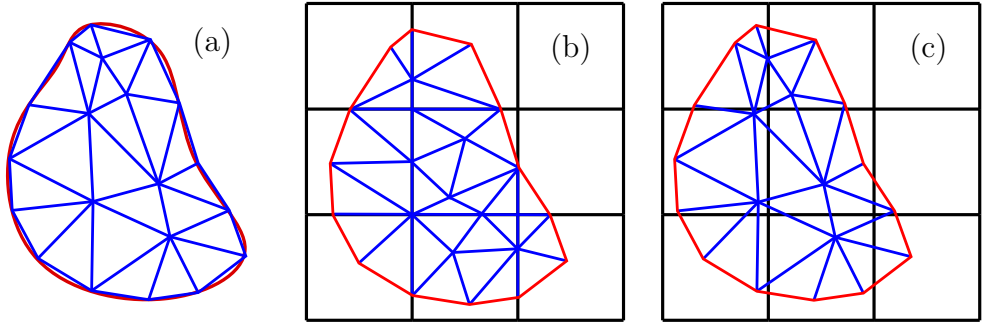

The finite element (FE) method \cite{Hughes87} is a widely-used approach for approximating solutions to partial differential equations (PDEs).  In summary, it approximates solutions to variational formulations of PDEs in function spaces defined piece-wise over a collection of disjoint elements, called a mesh, whose closures cover (an approximation of) the PDE domain.  Typically, each element is diffeomorphic to one of a few standard reference elements, like a unit right triangle in 2D or a tetrahedron in 3D, and the boundary of the (approximate) PDE domain coincides with a collection of element boundary parts, e.g., edges of elements in 2D and faces of elements in 3D. An example of such a ``boundary-fitted'' mesh is shown in Figure \ref{fig:immersed_meshes} (a).  However, the approximation power and regularity of FE function spaces depend on the sizes, shapes, and connectivity of the elements, placing nontrivial geometric and topological constraints on the mesh generation process.  These constraints are often quite challenging to satisfy, especially for the complicated PDE domains endemic to engineering analysis.  Thus, mesh generation remains a pervasive bottleneck in practical FE analysis \cite{Hardwick2005}.  

These meshing constraints become {\em much} easier to satisfy if one removes the requirement that the PDE boundary be a collection of element boundary parts, instead allowing it to cut arbitrarily through element interiors.  We refer to such meshes as ``unfitted'', and call numerical methods using unfitted meshes ``immersed boundary methods'' --- or simply ``immersed methods'' --- in an extrapolation of terminology originally applied to early work by Peskin \cite{Peskin72}.  For instance, it is easy to construct highly-regular function spaces with excellent approximation properties on structured tensor-product grids of elements, and it would be convenient to simply immerse an arbitrary complex PDE domain geometry into such a grid without needing its boundary to coincide with element boundary parts.

The earliest immersed methods amounted, essentially, to worsening the geometric approximation of the PDE domain, leading to low-order accuracy and the (still widespread) belief that immersed methods necessarily represent a trade-off between order of accuracy and ease of meshing.  However, over the past few decades, higher-order immersed methods have emerged based on the insight that one only needs an accurate representation of the domain geometry for numerical quadrature of variational forms from the FE method and not necessarily for constructing the approximate function spaces.  These quadrature-based immersed methods serve as our point of departure in the present work.

Quadrature-based immersed methods typically integrate weak forms by constructing a fitted mesh of the PDE domain and applying Gaussian quadrature within each element.  For instance, in Figure \ref{fig:immersed_meshes} (b), the fitted mesh of blue triangles would be used for quadrature, while the unfitted black structured mesh would be used for defining function spaces.  However, quadrature imposes {\em far} fewer geometric and topological constraints than construction of function spaces.  Further, even brute-force high-resolution/low-order quadrature meshes (such as those in early iterations of the finite cell method \cite{Parvizian07}) will not add degrees of freedom to the discrete function space and may remain efficient in applications whose cost is dominated by algebraic solution of the discretized PDE.  Thus, an easy-to-generate, ``foreground'' mesh for quadrature can be superposed on an unfitted ``background'' mesh for defining function spaces to obtain high-order-accurate immersed methods.  Some examples of accurate quadrature-based immersed methods from the literature include the finite cell approach mentioned earlier \cite{Schillinger2015} and various cut finite element (CutFEM) schemes \cite{Burman2015}.  One can also view the original material point method (MPM) \cite{Sulsky1995} as a quadrature-based immersed method using a low-order quadrature rule for the mass measure of a solid body. 

A historically-separate line of research aimed at circumventing the FE meshing problem is isogeometric analysis (IGA) \cite{Hughes05a,CoHuBa09} where function spaces of smooth splines from computer-aided design (CAD) are directly used to approximate PDE solutions.  However, these spline spaces often have {\em more} topological constraints on their construction than classical FE spaces, and industrial CAD models often represent complex geometries by ``trimming'' \cite{Marussig2018} unfitted rectangular meshes called patches to obtain the desired shapes.  Thus, IGA is hardly an alternative to immersed methods; on the contrary, it has only intensified interest in them!  The symbiotic union of isogeometric and immersed methods is referred to as ``immersogeometric analysis'' (IMGA) \cite{Kamensky2015} and the case of defining smooth spline spaces on a background mesh is sometimes also called CutIGA \cite{Elfverson2018} by analogy to CutFEM.  

An optimistic reading of the published results on immersed finite element and isogeometric methods might lead one to conclude that the FE mesh-generation problem is nearly solved.  However, the majority of work on immersed methods has relied on specialized research codes whose code architectures differ greatly from existing finite element or isogeometric solvers.  Thus, adoption of immersed methods in practice has been very limited.  A notable attempt at general-purpose implementation of quadrature-based immersed methods in an existing FE code was the MultiMesh \cite{Johansson2019} functionality developed within the FEniCS \cite{Logg2012} software framework.  However, the invasive nature of the implementation made it difficult to complete and maintain, and MultiMesh functionality was not retained in the FEniCSx redesign of the software.

A similar barrier to entry once existed for IGA, but a major advance in the design of IGA software was the concept of B\'{e}zier extraction \cite{BSEH11,SBVSH11} where smooth spline basis functions are represented as linear combinations of less-regular Bernstein basis functions.  This was later reformulated as {\em interpolation} of IGA basis functions using Lagrange FE bases, referred to as ``Lagrange extraction'' \cite{Schillinger2016} and allowing for direct reuse of FE software implementing Lagrange bases for IGA (e.g., \cite{Kamensky2019,Tirvaudey2020}).  {\bf The central idea of the present work is to build on the concept of Lagrange extraction by interpolating basis functions defined on an unfitted background mesh into a Lagrange FE space defined on a fitted foreground mesh where quadrature is performed.}  In this setting, the background basis does not need to be a smooth spline space; it could also be a classical FE space defined on an unstructured mesh.  This not only permits reuse of existing FE software, but also suggests some efficient approximations, leading to a novel class of immersed methods.

An earlier related idea is the projection-based adaptive Gaussian quadrature (PAGQ) implementation of the finite cell method \cite{Liu2020} where shape functions are defined on quadrature sub-cells.  However, we now push the concept further by considering the possibility of {\em approximate} Lagrange extraction.  While extraction-based IGA and PAGQ exactly reproduce basis functions through interpolation, we explore the possibility of allowing this interpolation to be approximate.  In situations like Figure \ref{fig:immersed_meshes} (b), where the foreground mesh is fitted to the intersection of each background element with the domain, it may be possible to capture all monomials of the background basis functions exactly using a sufficiently high-order foreground Lagrange space.  However, one may wish to reduce the order of the foreground space to reduce the cost of assembling algebraic equations.  Such an approximation is especially useful if the foreground elements are non-affine mappings of reference elements to capture the domain geometry with higher-order accuracy.  Further, one may consider more radical approximations, like Figure \ref{fig:immersed_meshes} (c), where the foreground mesh is unrelated to the background mesh beyond having a proportional quasi-uniform element size.  Notably, Gaussian quadrature on such a mesh would {\em not} be accurate for the background basis functions (which are non-smooth within foreground elements), but it {\em is} accurate for their interpolations into a foreground Lagrange function space.  Such an interpolation-based immersed method can be viewed as a generalization of the convected particle domain integration (CPDI) variant of the MPM \cite{Sadeghirad2011} proposed earlier for simulating large-strain solid mechanics, where our foreground elements are analogous to CPDI ``particles'' that carry their own shape functions.

The remainder of this article is structured as follows. Section \ref{sec:Notation} provides a concise list of notation and terminology for reference while reading subsequent sections (and may be skipped on an initial reading). Section \ref{sec:Immersed} introduces the concept of interpolation-based immersed methods using a scalar elliptic model PDE along with a sketch of how a basic error analysis of the method might proceed.  Section \ref{sec:Implementation} then describes how to implement interpolation-based immersed analysis for arbitrary PDE systems in a way that efficiently reuses existing FE software.  We then explore the accuracy of the method for a variety of PDEs through numerical experiments in Section \ref{sec:NumericalResults} and apply it to analysis of shell structures modeled geometrically by trimmed spline surfaces in Section \ref{sec:Application}.  Section \ref{sec:conclusion} draws conclusions and outlines future work.

\section{Notation and Terminology}
\label{sec:Notation}

The use of multiple meshes with different function spaces can lead to ambiguity in terms like ``element size'' and ``basis function'', so we must be careful to distinguish between background, foreground, and interpolated entities through precise terminology and notation.  We shall define notations more comprehensively as they arise in context later, but we also summarize them here for convenient reference:
\begin{itemize}
    \item $h$ : background mesh element size
    \item $\eta$ : foreground mesh element size
    \item $N_i$ : background basis functions
    \item $\widehat{N}_i$ : interpolated background basis functions
    \item $\phi_i$ : foreground basis functions
    \item $k$ : polynomial degree of the background FE space
    \item $\kappa$ : polynomial degree of the foreground FE space
    \item $\widehat{k}$ : limiting degree $\min\{k,\kappa\}$ of the interpolated background FE space
    \item $\mathcal{P}$ : interpolation into the foreground FE space
    \item $\mathcal{Q}$ : interpolation into the background FE space
    \item $\widehat{\mathcal{Q}}$ : interpolation $\mathcal{P}\circ\mathcal{Q}$ into the interpolated background FE space
\end{itemize}
We also define terms to distinguish the mesh immersions shown in Figures \ref{fig:immersed_meshes} (b) and \ref{fig:immersed_meshes} (c):
\begin{itemize}
    \item {\bf Background-fitted foreground mesh:} A foreground mesh which is fitted not just to the domain geometry's boundary but also to interior facets of the background mesh in the sense that the closure of the intersection of a background element with the PDE domain is the closure of a union of foreground elements, as shown in Figure \ref{fig:immersed_meshes} (b).
    \item {\bf Background-unfitted foreground mesh:} A foreground mesh such as the one shown in Figure \ref{fig:immersed_meshes} (c) whose elements are fitted to the domain geometry but not to interior facets of the background mesh.  
\end{itemize}
We may sometimes use the terms ``fitted'' and ``unfitted'' without specifying the prefix ``boundary-'' or ``background-'', but only when the omitted prefix should be clear from context.

\section{Immersed finite element analysis for a scalar elliptic model problem}
\label{sec:Immersed}

This section defines our proposed interpolation-based immersed method for a model scalar elliptic problem.  However, we will apply it to several more complicated problems in subsequent sections.  We define our model PDE in Section \ref{sec:model-pde}, then define a variational formulation suitable for immersed methods in Section \ref{sec:model-discrete}.  Section \ref{sec:model-quadrature-based} reviews the formulation and challenges of quadrature-based immersed methods, and Section \ref{sec:MainIdea} defines the construction of a function space used for interpolation-based immersed analysis. In Section \ref{sec:Approximation} we outline an {\em a priori} error estimate for the case of a background-fitted foreground mesh.

\subsection{The model problem}\label{sec:model-pde}

We choose as our model PDE the Poisson equation with Dirichlet boundary conditions.  Let $\Omega \subset \mathbb{R}^d$ be a spatial domain of interest.  For simplicity of exposition, we assume that $\Omega$ is polyhedral and can be decomposed into affine simplicial elements exactly.  The strong form of the problem is as follows: Find $u : \Omega \rightarrow \mathbb{R}$ such that
\begin{equation}
\begin{split}
\label{eq:poisson-strong}
    -\Delta u = f  ~~~ & \text{ in } \Omega \text{ ,} \\
            u = g  ~~~ & \text{ on } \partial \Omega \text{ .}
\end{split}
\end{equation}
The weak form we use as a starting point for formulating numerical methods is: Find $u \in \mathcal{V}_g$ such that

\begin{equation}
\int_{\Omega} \nabla u \cdot \nabla v d\Omega = \int_{\Omega} f v d\Omega\quad\forall v\in\mathcal{V}_0\text{ ,}\label{eq:weak-form-poisson}
\end{equation}
where
\begin{align}
\mathcal{V}_g &:= \left\{ v \in H^1(\Omega): v|_{\partial \Omega} = g \right\}, \\
\mathcal{V}_0 &:= \left\{ v \in H^1(\Omega): v|_{\partial \Omega} = 0 \right\},
\end{align}
in which restriction is understood in the sense of trace.

\subsection{Discretization of the model problem}\label{sec:model-discrete}
The most common discretization of the weak problem \eqref{eq:weak-form-poisson} is a conforming one employing finite dimensional subsets of $\mathcal{V}_g$ and $\mathcal{V}_0$.  However, the construction of these subsets typically relies on a boundary-fitted mesh to conform to the boundary conditions.  Many immersed methods therefore use nonconforming discretizations which enforce the Dirichlet boundary condition in a weak sense.  In particular, we employ Nitsche's method \cite{Nitsche71} for weak enforcement of Dirichlet boundary conditions, leading to the discrete problem:
Find $u^h \in \mathcal{V}^h$ such that, $\forall v^h\in\mathcal{V}^h$,
\begin{align}
\int_{\Omega} \nabla u^h \cdot \nabla v^h d\Omega - \int_{\partial \Omega} \left( \nabla u^h \cdot \bm{n}\right) v^h d\Gamma \mp \int_{\partial \Omega} \left( \nabla v^h \cdot \bm{n}\right) u^h d\Gamma + \int_{\partial \Omega} \frac{C_\text{pen}}{h} u^h v^h d\Gamma \nonumber \\
= \int_{\Omega} f v^h d\Omega \mp \int_{\partial \Omega} \left( \nabla v^h \cdot \bm{n}\right) g d\Gamma + \int_{\partial \Omega} \frac{C_\text{pen}}{h} g v^h d\Gamma\text{ ,}\label{eq:poisson-disc-generic}
\end{align}
where the discrete space $\mathcal{V}^h$ is a finite-dimensional subspace of $H^1(\Omega)$ with no extra restrictions on $\partial\Omega$, $h$ is a length scale indicating the refinement level of the discretizations (e.g., the element diameter in a standard boundary-fitted discretization), and $C_\text{pen} \geq 0$ is a dimensionless constant.  The ``$\mp$'' in \eqref{eq:poisson-disc-generic} must be selected consistently on both sides of the equation, and toggles between the symmetric (``$-$'') and non-symmetric (``$+$'') variants of the Nitsche method.  The original symmetric variant proposed by \cite{Nitsche71} remains well-posed for $C_\text{pen}$ sufficiently large, with the lower bound emanating from trace-inverse estimates (but often chosen on an {\em ad hoc} basis in practice).  The non-symmetric variant is well-posed for $C_\text{pen}$ arbitrarily small (or zero), which has led some authors to advocate its use in immersed discretizations \cite{Schillinger2016b} to sidestep technical concerns over the appropriate values for $C_\text{pen}$ and $h$.  The non-symmetric variant retains optimal accuracy in the $H^1$ norm, which is often more relevant to applications like stress analysis.  However, the symmetric variant is optimal in both $H^1$ and $L^2$ norms, and is currently more widely used in practice.

\subsection{Quadrature-based immersed discretization}\label{sec:model-quadrature-based}
In quadrature-based immersed methods, such as those reviewed in Section \ref{sec:introduction}, the discrete space $\mathcal{V}^h$ is defined by simply restricting functions in an FE (or isogeometric) space defined on a background mesh to the PDE domain $\Omega$:
\begin{equation}
    \mathcal{V}^h = \operatorname{span}\left\{N_i\vert_{\Omega}\right\}_{i=1}^n\text{ ,}
\end{equation}
where $\{N_i\}_{i=1}^n$ are the basis functions of the background-mesh FE space whose supports intersect $\Omega$.  The principal challenge in implementing such a method is to accurately compute integrals of the forms
\begin{equation}\label{eq:split-integrals}
    \int_\Omega(\cdots)\,d\Omega = \sum_{E=1}^{n_{\text{el}}}\int_{\Omega_E\cap\Omega}(\cdots)\,d\Omega\quad\text{and}\quad\int_{\partial\Omega}(\cdots)\,d\Gamma = \sum_{E=1}^{n_\text{el}}\int_{\partial\Omega\cap\overline{\Omega_E}}(\cdots)\,d\Gamma\text{ ,}
\end{equation}
where $\Omega_E$ is the $E^\text{th}$ element of the background mesh and $n_\text{el}$ is the total number of background mesh elements.  One can use standard Gaussian quadrature on $\Omega_E$ when $\Omega_E\subset\Omega$ or on boundary parts of $\Omega_E$ when the background mesh is fitted to $\partial\Omega$.  Various schemes have been proposed for integrating nontrivial intersections, but we focus in this paper on those which decompose $\Omega_E\cap\Omega$ into a collection of quadrature elements, as shown in Figure \ref{fig:immersed_meshes} (b), and we refer to the mesh formed by these quadrature elements as a background-fitted foreground mesh.  With such a mesh, one can then decompose integrals like
\begin{equation}
    \int_\Omega(\cdots)\,d\Omega = \sum_{e=1}^{\nu_\text{el}}\int_{\omega_e}(\cdots)\,d\Omega\quad\text{and}\quad\int_{\partial\Omega}(\cdots)\,d\Gamma = \sum_{e=1}^{\nu_\text{el}}\int_{\partial\Omega\cap\partial\omega_e}(\cdots)\,d\Gamma\text{ ,}
\end{equation}
where $\{\omega_e\}_{e=1}^{\nu_\text{el}}$ are the elements of the foreground mesh and standard Gaussian quadrature rules can be accurately used on each one of its boundary parts.

While mathematically-appealing, this presents a major challenge from a software development standpoint as now the mesh over which quadrature rules are defined is a separate data structure from the mesh on which the basis functions are constructed.  Implementation of such methods therefore requires either custom software to be written from scratch or complicated invasive modifications to be made in mature FE or IGA software designed for standard boundary-fitted analysis.  The primary motivation of the present work is to circumvent these difficulties while retaining the high-order accuracy of quadrature-based immersed methods. 

\subsection{Interpolation-based immersed discretization}
\label{sec:MainIdea}
The main idea of the present work is to use a discrete space that consists of interpolations of background-mesh basis functions on a Lagrange FE space defined over the foreground mesh:
\begin{equation}
    \mathcal{V}^h = \operatorname{span}\left\{\widehat{N}_i\right\}_{i=1}^n\text{ ,}
\end{equation}
where
\begin{equation}
    \widehat{N}_i = \mathcal{P}N_i := \sum_{j=1}^\nu N_i(\boldsymbol{x}_j)\phi_j\text{ ,}
\end{equation}
in which $\mathcal{P}$ is the operator that interpolates a function into the Lagrange FE space with basis $\{\phi_j\}_{j=1}^\nu$ and $\boldsymbol{x}_j$ is the nodal point associated with foreground Lagrange basis function $\phi_j$, i.e.,
\begin{equation}
    \phi_i(\boldsymbol{x}_j) = \delta_{ij}\text{ .}
\end{equation}
The space $\mathcal{V}^h$ is a subset of the Lagrange FE space defined on the foreground mesh, so its functions can be accurately integrated using Gaussian quadrature rules on each foreground element.  However, unlike in quadrature-based immersed methods, this remains true even with a background-unfitted foreground mesh, like the one shown in Figure \ref{fig:immersed_meshes} (c), where, for an arbitrary background element $\Omega_E$, $\overline{\Omega_E\cap\Omega}$ is not necessarily equal to the closure of a union of foreground elements.

For either background-fitted or background-unfitted foreground meshes, the basis $\{\widehat{N}_i\}_{i=1}^{n}$ exhibits the following properties:
\begin{enumerate}
\item If  $\{N_i\}_{i=1}^{n}$ forms a partition of unity, then $\{\widehat{N}_i\}_{i=1}^{n}$ does as well. 
\item If the background function space spans polynomials of degree $k$ and the foreground function space spans polynomials of degree $\kappa$, then $\{\widehat{N}_i\}_{i=1}^{n}$ spans polynomials of degree $\widehat{k} :=\min\{k,\kappa\}$.
\item If both $\{N_i\}_{i=1}^{n}$ and $\{\phi_i\}_{i=1}^{^\nu}$ have local support, then $\{\widehat{N}_i\}_{i=1}^{n}$ also has local support. 
\item The continuity of $\{\widehat{N}_i\}_{i=1}^{n}$ is at least that of $\{\phi_i\}_{i=1}^{^\nu}$.  Note that this may be lower than the continuity of $\{N_i\}_{i=1}^{n}$.  
\end{enumerate}
Moreover, for the specific case of a background-fitted foreground mesh, the interpolated background basis $\{\widehat{N}_i\}_{i=1}^{n}$ is equivalent to the original background basis $\{N_i\}_{i=1}^{n}$ for a sufficiently high foreground polynomial degree $\kappa$.  If the background function space consists of piecewise polynomials of degree $k$ defined over a simplicial background mesh, the original background basis is recovered with a Lagrange interpolation of degree $\kappa \geq k$ onto a simplicial background-fitted foreground mesh.  If the background function space consists of piecewise tensor-product polynomials of degree $k$ defined over a tensor-product background mesh, the original background basis is recovered with a Lagrange interpolation of degree $\kappa \geq k^d$ onto a simplicial background-fitted foreground mesh.

\afterpage{

\begin{figure}
    \centering
     \hspace{15pt}
	 \begin{subfigure}[b]{0.315\linewidth}
	 \centering
    \includegraphics[width=\linewidth]{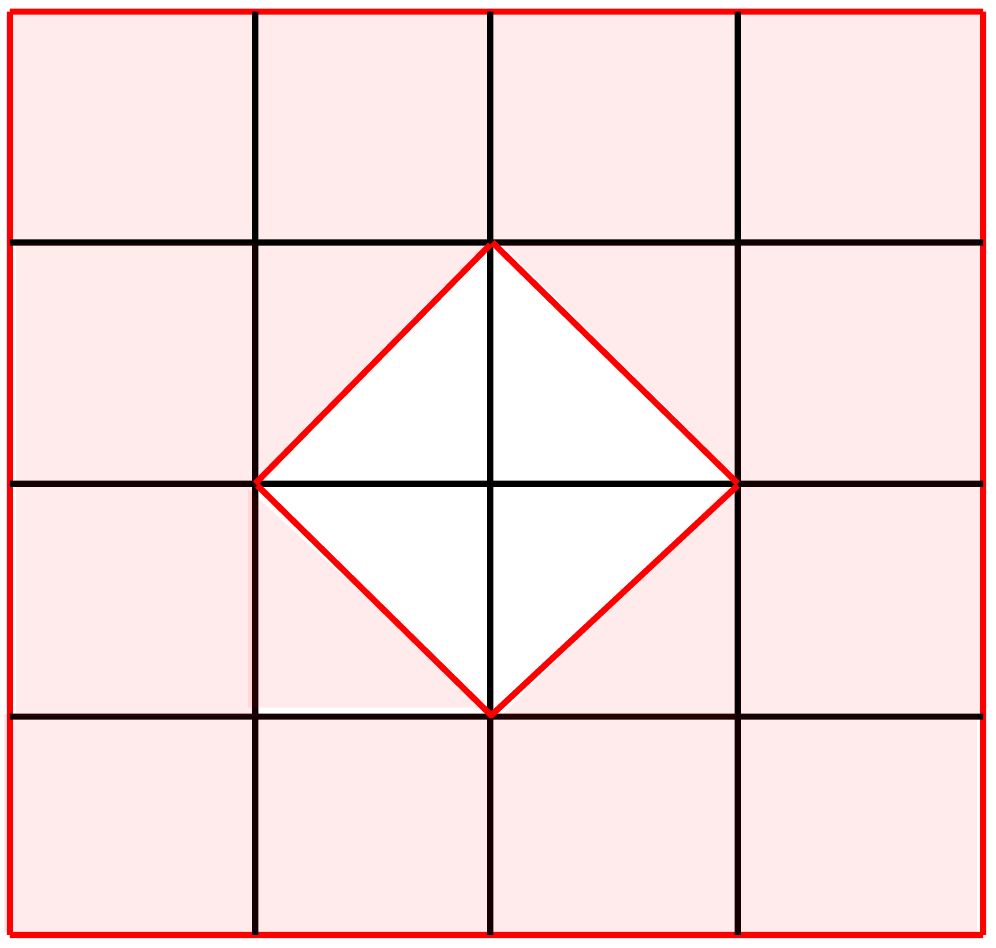} (a)
    \label{fig:basis_function_a}
  \end{subfigure}
    \hspace{45pt}
	 \begin{subfigure}[b]{0.40\linewidth}
	 \centering
    \includegraphics[width=\linewidth]{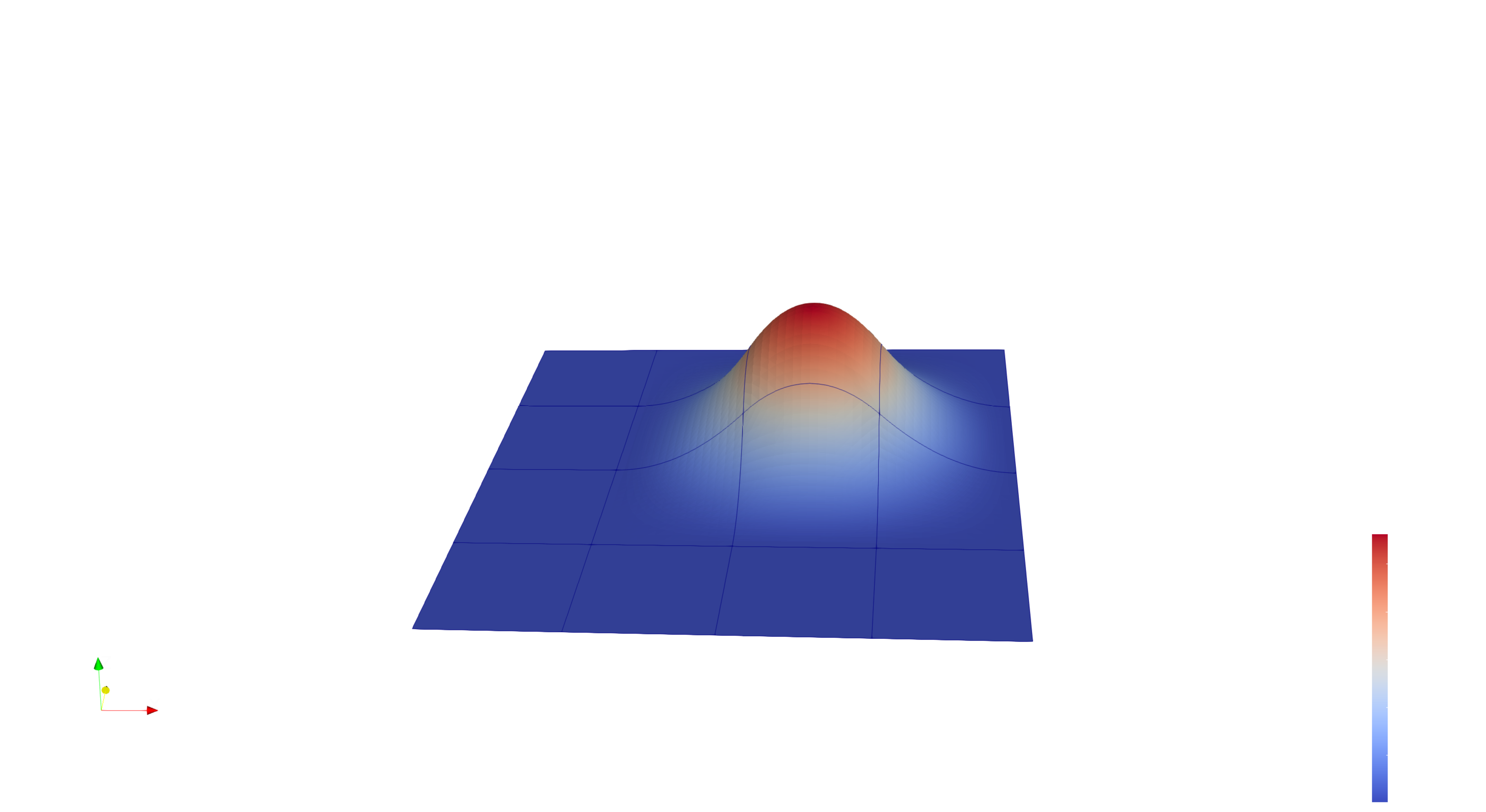} (b)
    \label{fig:basis_function_b}
  \end{subfigure}
  \centering
	 \begin{subfigure}[b]{0.33\linewidth}
	 \centering
    \includegraphics[width=\linewidth]{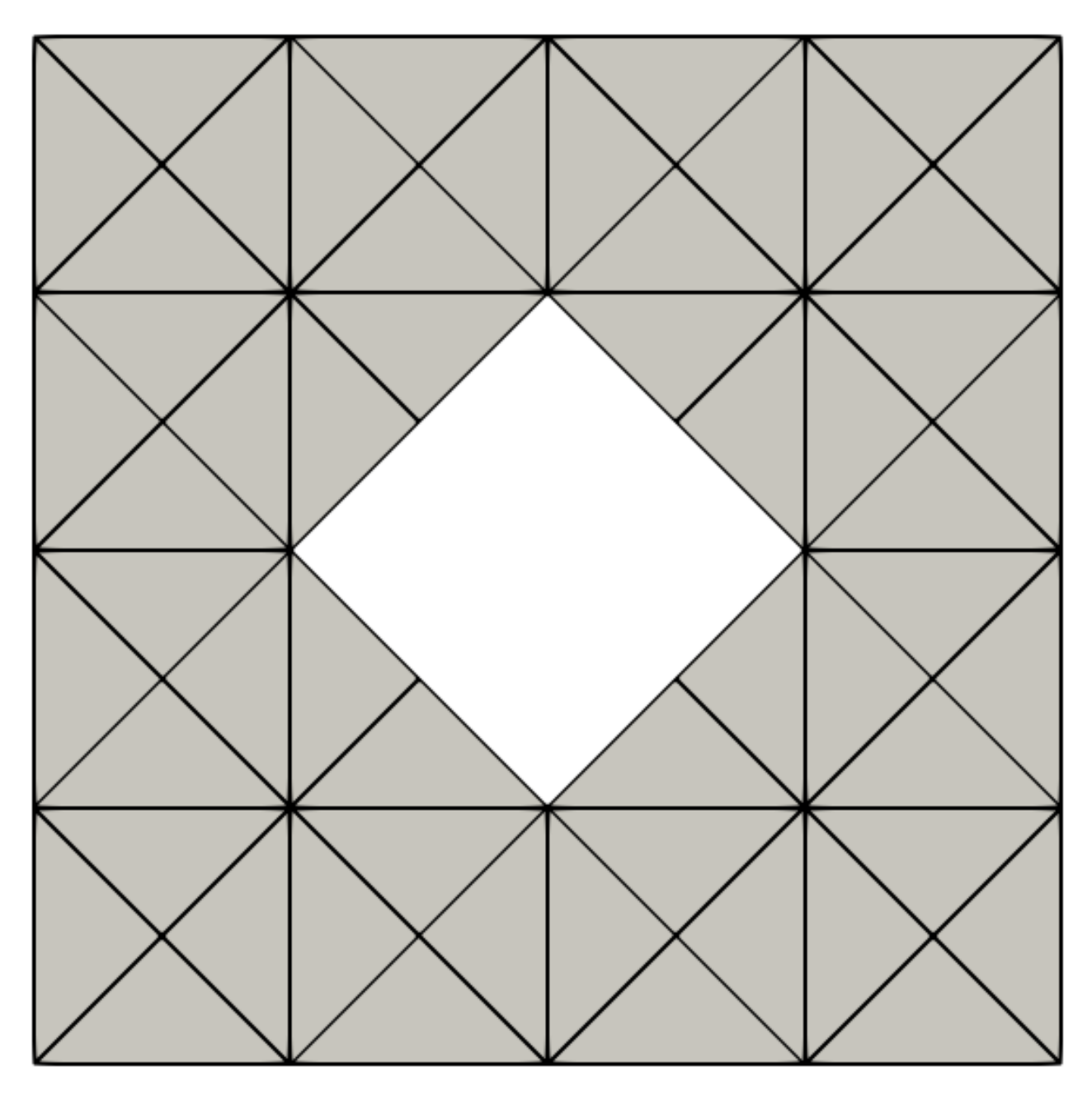} (c)
    \label{fig:basis_function_c}
  \end{subfigure}
    \hspace{59pt}
	 \begin{subfigure}[b]{0.33\linewidth}
	 \centering
    \includegraphics[width=\linewidth]{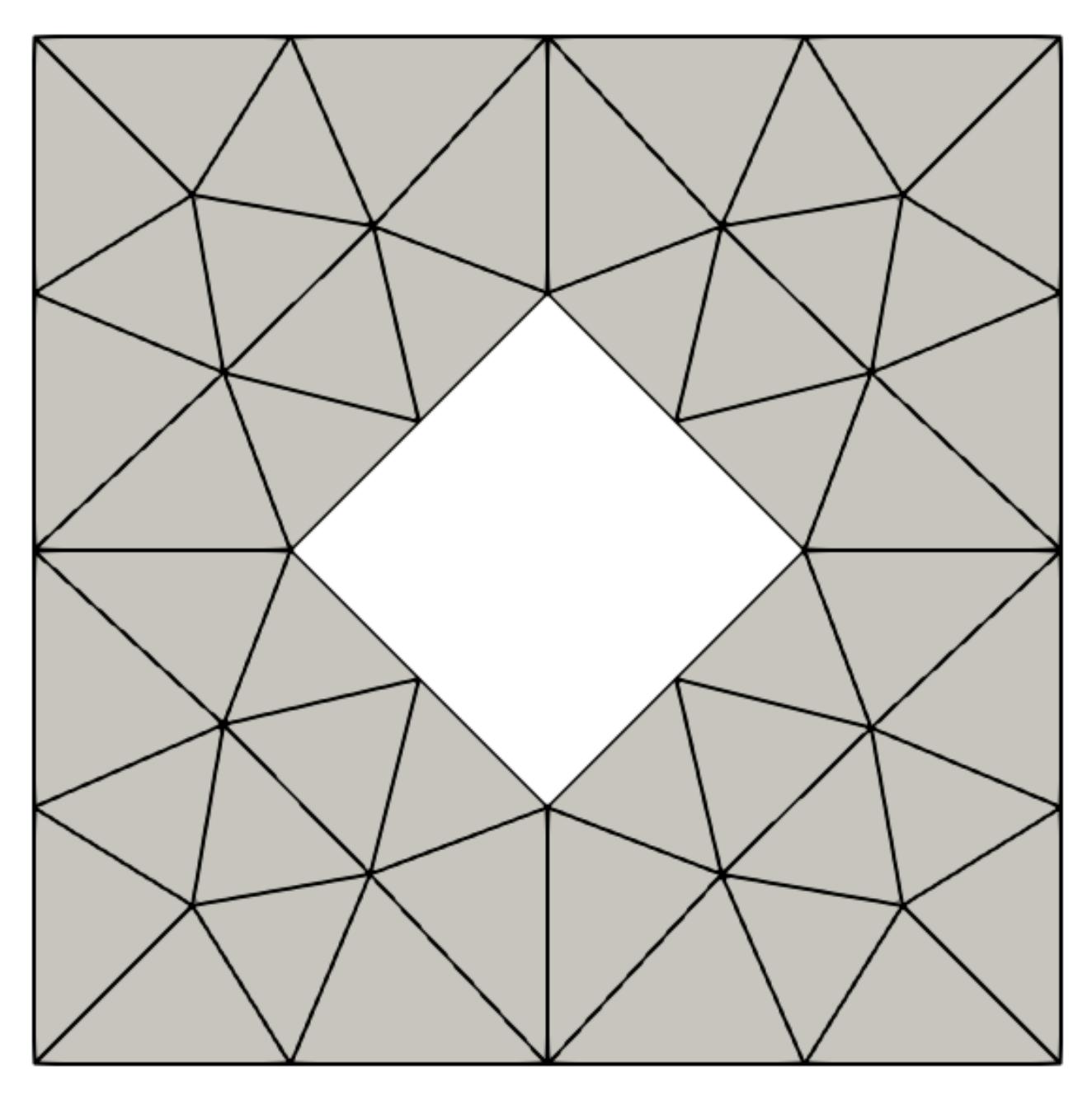} (d)
    \label{fig:basis_function_d}
  \end{subfigure}
    \begin{subfigure}[b]{0.40\linewidth}
    \centering
    \includegraphics[width=\linewidth]{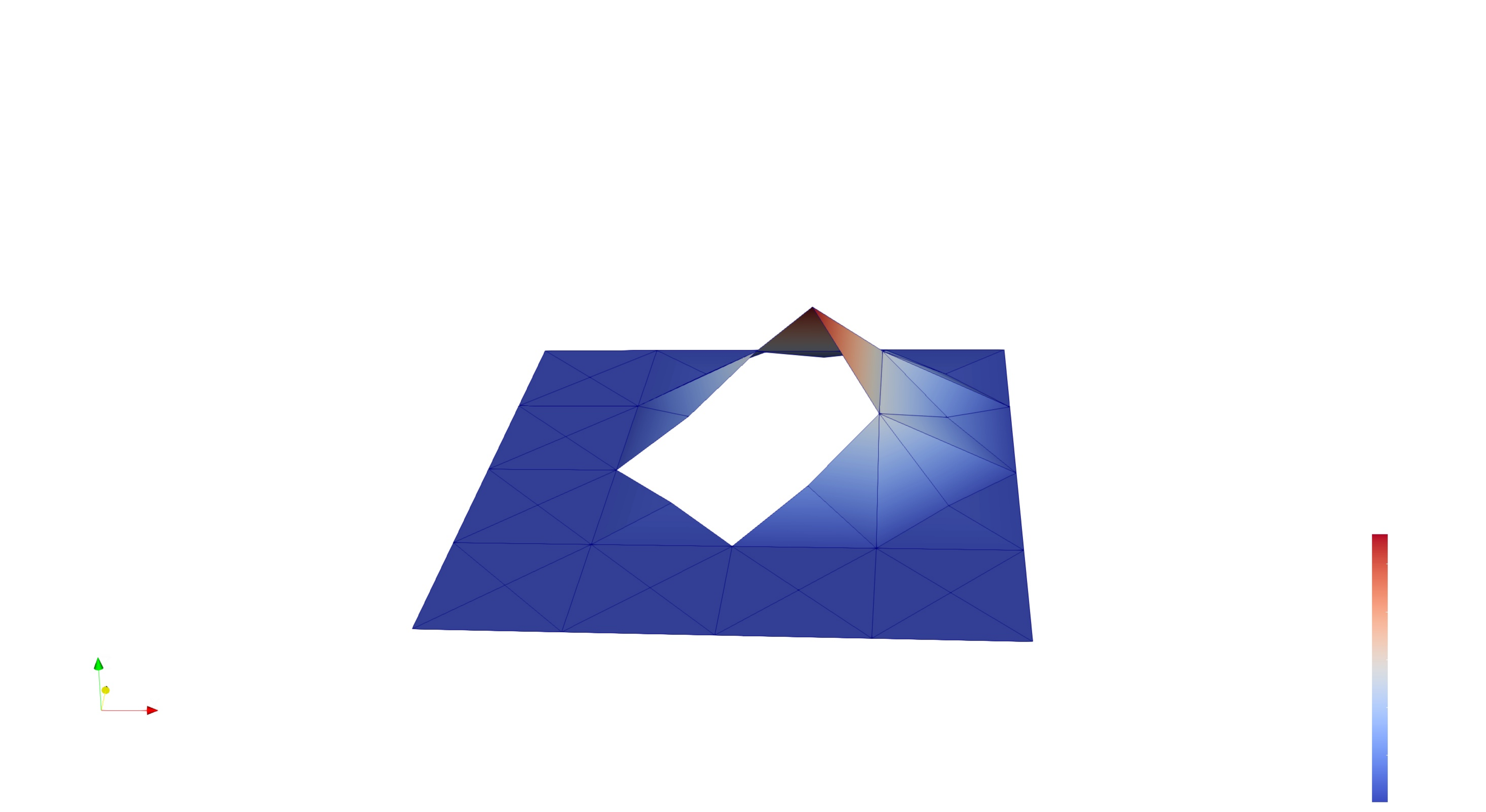} (e)
    \label{fig:basis_function_e}
  \end{subfigure}
    \hspace{25pt}
	 \begin{subfigure}[b]{0.40\linewidth}
	 \centering
    \includegraphics[width=\linewidth]{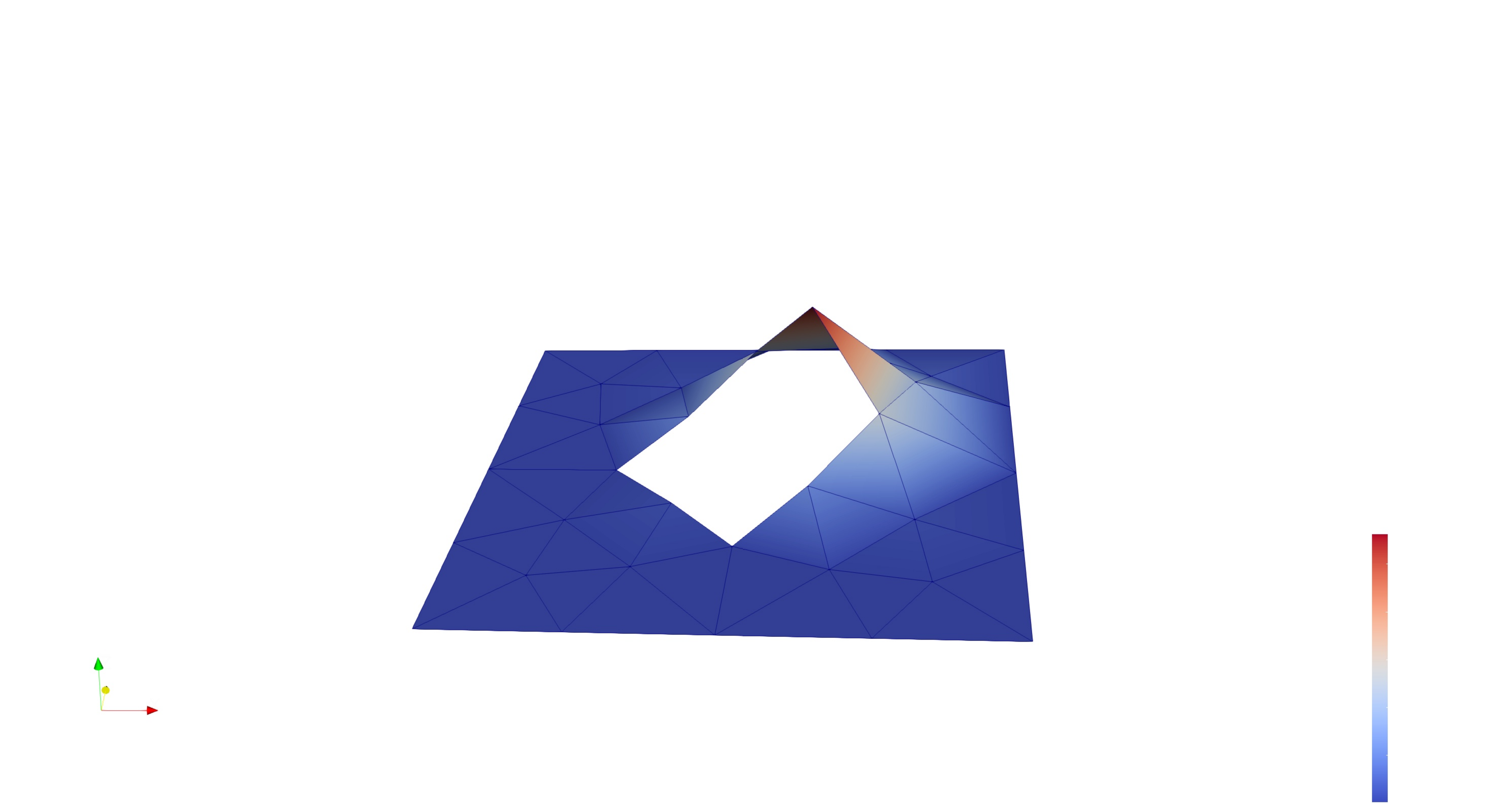} (f)
    \label{fig:basis_function_f}
  \end{subfigure}
  \begin{subfigure}[b]{0.40\linewidth}
	 \centering
    \centering
    \includegraphics[width=\linewidth]{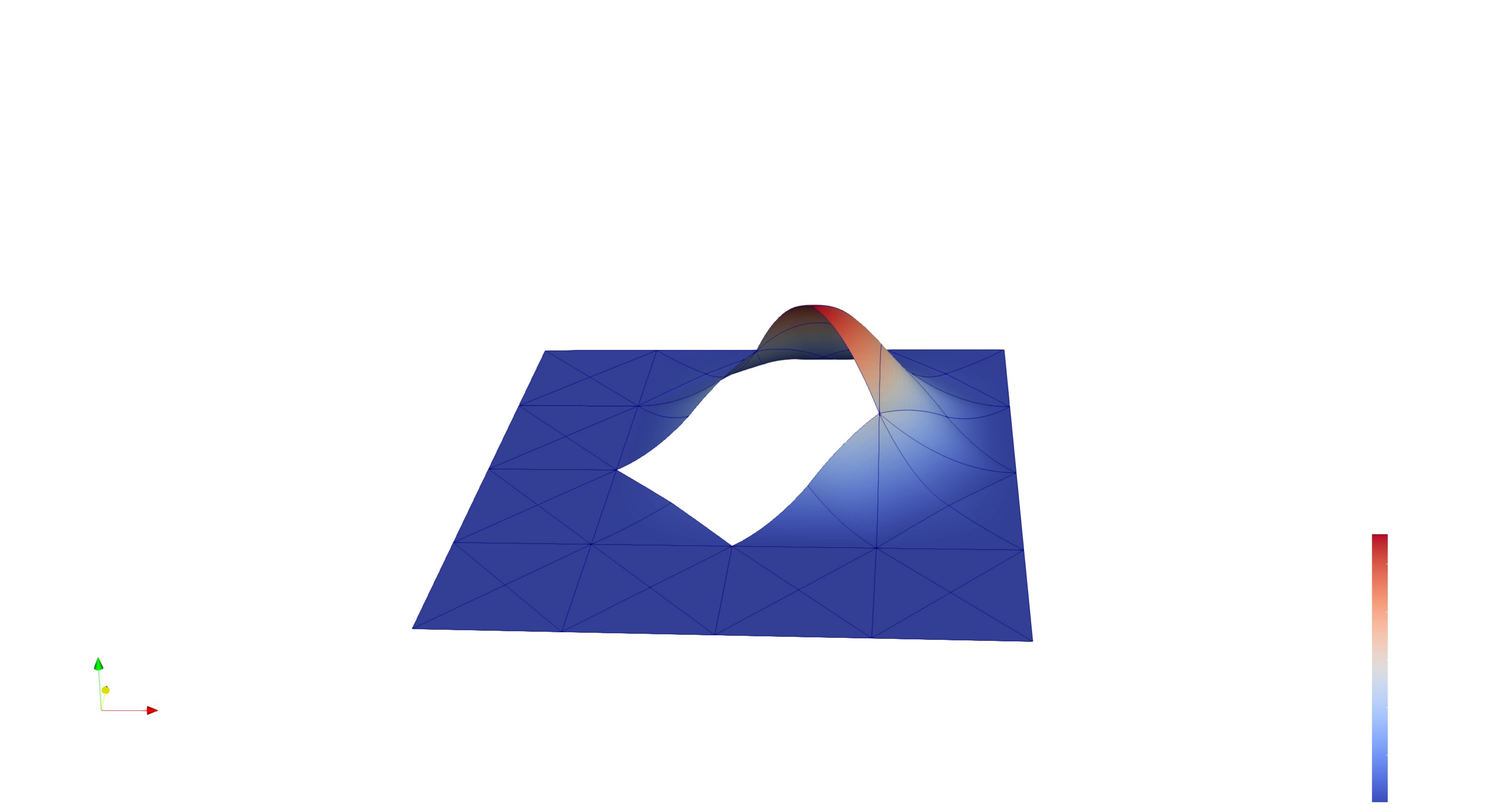} (g)
    \label{fig:basis_function_g}
  \end{subfigure}
    \hspace{25pt}
	 \begin{subfigure}[b]{0.40\linewidth}
	 \centering
    \includegraphics[width=\linewidth]{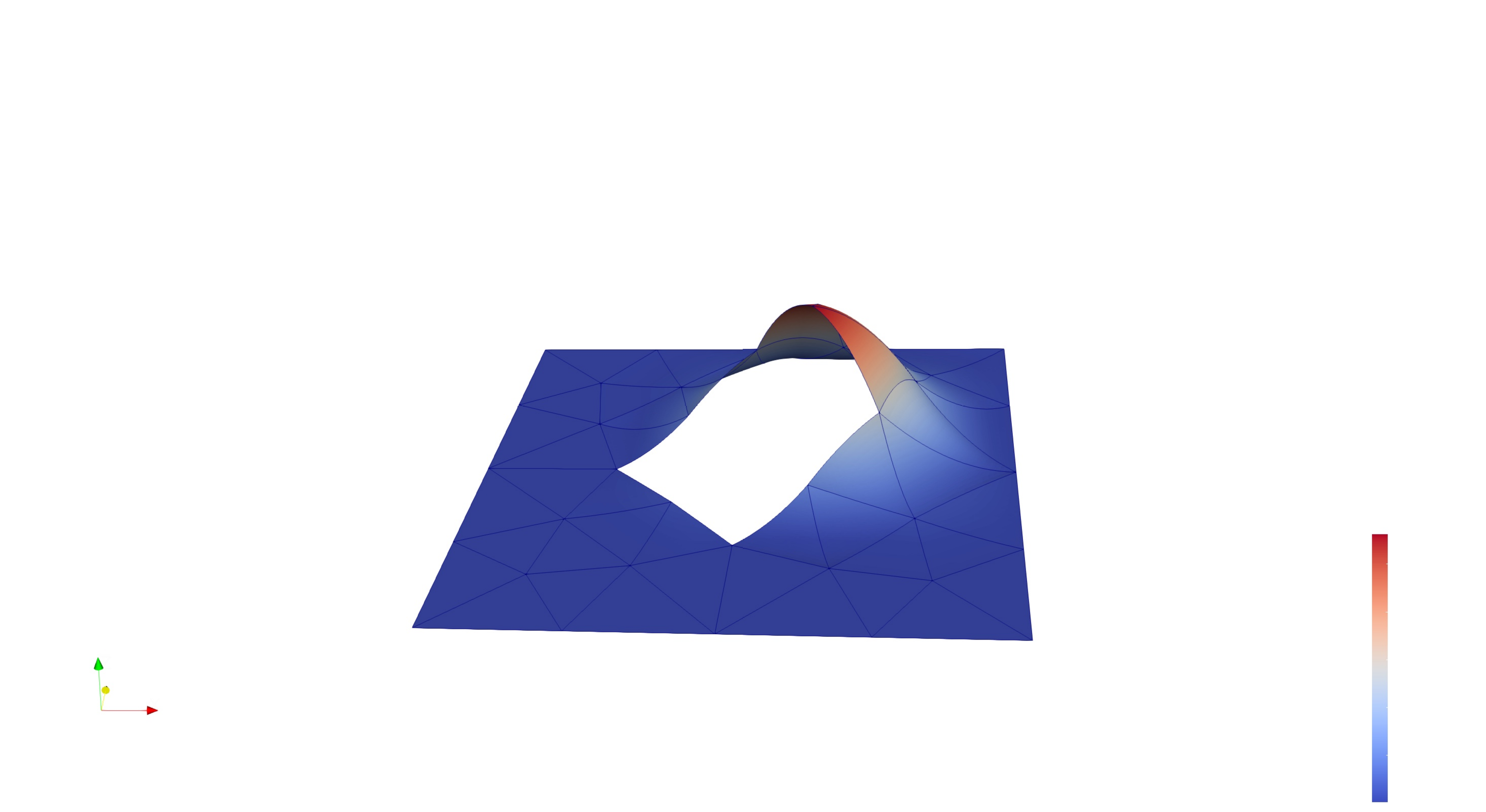} (h)
    \label{fig:basis_function_h}
  \end{subfigure}
    \caption{(a) Geometric domain (red) and background mesh (black).  (b) Background basis function.  (c) Background-fitted foreground mesh.  (d) Background-unfitted foreground mesh. (e) Linear Lagrange interpolation of background basis function on background-fitted mesh.  (f) Linear Lagrange interpolation of background basis function on background-unfitted mesh.  (g) Quadratic Lagrange interpolation of background basis function on background-fitted mesh.  (h) Quadratic Lagrange interpolation of background basis function on background-unfitted mesh.}
    \label{fig:basis_function}
\end{figure}

\clearpage
}

To better convey the concept of interpolated background basis function, several interpolations of an example background basis function are displayed in Figure \ref{fig:basis_function}.  Figure \ref{fig:basis_function} (a) shows the geometric domain and background mesh in red and black respectively, while Figure \ref{fig:basis_function} (b) shows the example background basis function.  The background basis function is a bi-quadratic ($k = 2$) B-spline basis function of maximal $C^{k-1}$ continuity defined over the background mesh.  Figure \ref{fig:basis_function} (c) and Figure \ref{fig:basis_function} (d) show representative background-fitted and background-unfitted foreground meshes.  Figure \ref{fig:basis_function} (e) and Figure \ref{fig:basis_function} (f) show linear ($\kappa = 1$) Lagrange interpolations of the example background basis function onto the background-fitted and background-unfitted foreground meshes.  These two interpolated background basis functions share the same shape as the background basis function, but they are both inaccurate approximations of the background basis function.  Alternatively, Figure \ref{fig:basis_function} (g) and Figure \ref{fig:basis_function} (h) show quadratic ($\kappa = 2$) Lagrange interpolations of the background basis function onto the background-fitted and background-unfitted foreground meshes.  These two interpolated background basis functions are much more accurate approximations of the background basis function than the interpolated background basis functions displayed in Figure \ref{fig:basis_function} (e) and Figure \ref{fig:basis_function} (f).

\subsection{Best approximation estimates for the background-fitted case}\label{sec:Approximation}
Having presented our new paradigm for immersed finite element and isogeometric analysis, we now derive estimates for the $H^1$-norm best approximation error
\begin{equation}
\inf_{v^h \in \operatorname{span}\left\{\widehat{N}_i\right\}_{i=1}^n} \left\Vert u - v^h\right\Vert_{H^1(\Omega)} \label{eq:bae_etimate}
\end{equation}
associated with the interpolated background FE space where $u$ is the exact solution of our model problem.  Best approximation estimates of the form given by \eqref{eq:bae_etimate} can be combined with coercivity and continuity results to arrive at \textit{a priori} error estimates for a given immersed finite element method (as done, e.g., for quadrature-based CutFEM with ghost penalty stabilization in \cite{Burman2015}).  In our mathematical analysis, we must make two assumptions regarding the foreground mesh:

\medskip
\noindent \textbf{Assumption 1:} \textit{The foreground mesh is background-fitted.}\\

\noindent \textbf{Assumption 2:} \textit{The foreground mesh is free of hanging nodes.}\\

\noindent Assumption 1 ensures that the background basis functions are smooth within every element of the foreground mesh, while Assumption 2 ensures that the interpolated background FE space is $H^1$-conforming and thus the $H^1$-norm best approximation error is well-defined.  Note that an immediate consequence of Assumption 1 is that the foreground mesh size $\eta$ is less than or equal to the background mesh size $h$.  We must make two further assumptions regarding the exact solution $u$ and the background FE space:\\

\noindent \textbf{Assumption 3:} \textit{The exact solution $u$ lies in the space $H^{\widehat{k}+1}(\Omega)$.}\\

\noindent \textbf{Assumption 4:} \textit{There exists an interpolation operator $\mathcal{Q}: H^1(\Omega) \rightarrow \operatorname{span}\left\{N_i\right\}_{i=1}^n$ onto the background FE space such that
\begin{equation}
\left\Vert v - \mathcal{Q} v \right\Vert_{H^1(\Omega)} \leq C_1 h^{\widehat{k}} \left\Vert v \right\Vert_{H^{\widehat{k}+1}(\Omega)} \label{eq:bg_interpolation_assumption}
\end{equation}
and
\begin{equation}
\sum_{e=1}^{\nu_\text{el}} \left\Vert \mathcal{Q} v \right\Vert^2_{H^{\widehat{k}+1}(\omega^e)} \leq C_2 \left\Vert v \right\Vert^2_{H^{\widehat{k}+1}(\Omega)} \label{eq:bg_continuity_assumption}
\end{equation}
for all $v \in H^{\widehat{k}+1}(\Omega)$ where $C_1 > 0$ and $C_2 > 0$ depend on the background FE space polynomial degree $k$ and the shape regularity of the background mesh but not the background mesh element size $h$.}\\

\noindent Assumption 3 holds provided the problem geometry and parameters are sufficiently smooth, while Assumption 4 holds when the background FE space consists of classical finite element functions or smooth spline functions provided there exists a continuous extension operator $E: H^{k+1}(\Omega) \rightarrow H^{k+1}(\mathbb{R}^d)$ (cf. \cite[(9)]{Burman2015}).

Now let us define $\widetilde{\mathcal{Q}} = \mathcal{P} \circ \mathcal{Q}$.  Then we can bound the $H^1$-norm best approximation error associated with the interpolated background FE space by the $H^1$-norm error associated with the interpolation operator $\widetilde{\mathcal{Q}}$:
\begin{equation}
\inf_{v^h \in \operatorname{span}\left\{\widehat{N}_i\right\}_{i=1}^n} \left\Vert u - v^h\right\Vert_{H^1(\Omega)} \leq \left\Vert u - \widetilde{\mathcal{Q}} u \right\Vert_{H^1(\Omega)}. \label{eq:bae_bound}
\end{equation}
By the triangle inequality,
\begin{equation}
\left\| u - \widetilde{\mathcal{Q}} u \right\|_{H^1(\Omega)} \leq \| u - \mathcal{Q} u \|_{H^1(\Omega)} + \| \mathcal{Q} u - \mathcal{P} \mathcal{Q} u \|_{H^1(\Omega)}\text{ .} \label{eq:bae_term_1}
\end{equation}
Invoking Assumption 3, we can immediately bound the first term on the right hand side of the above equation:
\begin{equation}
\| u - \mathcal{Q} u \|_{H^1(\Omega)} \leq C_1 h^{\widehat{k}} \left\Vert u \right\Vert_{H^{{\widehat{k}}+1}(\Omega)}. \label{eq:bae_term_2}
\end{equation}
To bound the second term, we first recognize that
\begin{equation}
\| \mathcal{Q} u - \mathcal{P} \mathcal{Q} u \|^2_{H^1(\Omega)} = \sum_{e=1}^{\nu_\text{el}} \| \mathcal{Q} u - \mathcal{P} \mathcal{Q} u \|^2_{H^1(\omega^e)}. \label{eq:bae_term_3}
\end{equation}
By Assumption 1, we have that $\mathcal{Q} u|_{\omega^e}$ is smooth for every foreground mesh element $\omega^e$.  Thus we have that
\begin{equation}
\| \mathcal{Q} u - \mathcal{P} \mathcal{Q} u \|^2_{H^1(\omega^e)} = C_3 \eta^{2\widehat{k}} \| \mathcal{Q} u \|^2_{H^{{\widehat{k}}+1}(\omega^e)} \label{eq:bae_term_4}
\end{equation}
for every foreground element $\omega^e$ where $C_3$ is a constant that depends on the foreground FE space polynomial degree $\kappa$ and the shape regularity of the foreground mesh but not the foreground mesh element size $\eta$.  We can immediately combine \eqref{eq:bae_term_3}, \eqref{eq:bae_term_4}, Assumption 1, and Assumption 4 to arrive at
\begin{equation}
\| \mathcal{Q} u - \mathcal{P} \mathcal{Q} u \|_{H^1(\Omega)} =\sqrt{ C_2 C_3 } h^{\widehat{k}} \| u \|_{H^{{\widehat{k}}+1}(\Omega)} \label{eq:bae_term_5}.
\end{equation}
Combining \eqref{eq:bae_bound}, \eqref{eq:bae_term_1}, \eqref{eq:bae_term_2}, and \eqref{eq:bae_term_3}, we finally arrive at the following estimate for the $H^1$-norm best approximation error associated with the interpolated background FE space:
\begin{equation}
\inf_{v^h \in \operatorname{span}\left\{\widehat{N}_i\right\}_{i=1}^n} \left\Vert u - v^h\right\Vert_{H^1(\Omega)} \leq C_{\textup{interp}} h^{\widehat{k}} \| u \|_{H^{{\widehat{k}}+1}(\Omega)} \label{eq:final_bae_estimate}
\end{equation}
where $C_{\textup{interp}} = C_1 + \sqrt{C_2 C_3}$.  It is immediately evident that the above estimate is optimal with respect to both the background mesh size $h$ and the limiting degree $\hat{k}$ of the interpolated background space.  
Note that the convergence rate does not improve by decreasing the foreground mesh size $\eta$ faster than the background mesh size $h$
or increasing the foreground polynomial degree $\kappa$ above that of the background polynomial degree $k$.\footnote{If there is geometry error in the foreground mesh, then there can be some benefit to over-refining $\eta$ near boundaries (cf. Section \ref{subsec:linearElasticity}), but geometry error is outside the scope of the analysis in this section.}

In the above analysis, we made the rather stringent assumption that the foreground mesh be background-fitted.  However, our later numerical experiments suggest that best approximation estimates of the form given by \eqref{eq:final_bae_estimate} also hold when the foreground mesh is background-unfitted.  This is perhaps not too unexpected as the interpolated background FE space is complete up to polynomial degree $\hat{k}$ even when the foreground mesh is background-unfitted, but our analysis does not extend to this case since we require that the background basis functions be smooth over each element of the foreground mesh.  We also made the assumption in our analysis that the foreground mesh be free of hanging nodes.  This was to ensure that the interpolated background FE space was $H^1$-conforming.  That being said, our later numerical experiments suggest we may not need an interpolated background FE space that is $H^1$-conforming provided that the original background FE space is $H^1$-conforming.  In particular, we are able to attain optimal convergence rates for the biharmonic problem using an interpolated background FE space that is only $H^1$-conforming when the original background FE space is $H^2$-conforming.  Finally, the best approximation estimate given by \eqref{eq:final_bae_estimate} depends on the shape regularity of the foreground mesh through the interpolation constant $C_{\text{interp}}$, but our later numerical experiments suggest that high quality results are attained with interpolation-based immersed finite element and isogeometric analysis even when the foreground mesh is of poor quality.

\section{Reusing FE software for interpolation-based immersed methods}
\label{sec:Implementation}

As each basis function $\widehat{N}_i$ of $\mathcal{V}^h$ in our proposed method is a linear combination of basis functions defined on the foreground mesh used for element-by-element assembly, one may implement interpolation-based immersed methods using existing FE software in a similar manner as one would implement IGA via Lagrange extraction \cite{Schillinger2016}.  The cited paper emphasizes an implementation based on modifying element-level shape function routines which could also be followed for immersed methods.  However, the method we shall outline here is based on global linear algebra operations applied to the full, assembled system of equations, similar to what is done in \cite{Kamensky2019,Tirvaudey2020} for IGA.  The latter approach requires fewer modifications to  the existing FE code and is referred to as ``non-invasive'' by \cite{Tirvaudey2020}.

\subsection{Non-invasive implementation using extraction matrices}
\label{subsec:exMat}

Consider a scalar variational problem of the form: Find $u^h\in\mathcal{V}^h$ such that, $\forall v^h\in\mathcal{V}^h$,
\begin{equation}
    a\left(u^h,v^h\right) = L\left(v^h\right)\text{ ,}
\end{equation}
where $a$ and $L$ are a bilinear and linear form, i.e., the left- and right-hand sides of the Nitsche formulation \eqref{eq:poisson-disc-generic} for the Poisson model problem. $\mathcal{V}^h = \operatorname{span}\{\widehat{N}_i\}$ is the solution space for the interpolation-based immersed method. The state variable field is discretized by the interpolated basis functions as
\begin{equation}
    u^h = \sum_{i=1}^nd_i\widehat{N}_i\text{ ,}
\end{equation}
where $\{d_i\}_{i=1}^n$ are unknown coefficients. Then we need to solve the linear system
\begin{equation}
    \mathbf{K}\mathbf{d} = \mathbf{F}\text{ ,}\label{eq:linear-system}
\end{equation}
where the global stiffness matrix and forcing vector are
\begin{equation}
    K_{ij} = a\left(\widehat{N}_j,\widehat{N}_i\right)\quad\text{and}\quad F_i = L\left(\widehat{N}_i\right)\text{ .}\label{eq:defn-k-f}
\end{equation}
In the interpolation-based immersed method,
\begin{equation}
    \widehat{N}_j = \sum_{i=1}^{\nu}\exMat_{ij}\phi_i\quad\text{with}\quad \exMat_{ij} = N_j(\boldsymbol{x}_i)\text{ .}
\end{equation}
$\exMat_{ij}$ projects the background functions onto the space of the foreground elements and can be obtained by evaluating the background basis ${N_j}$ at the nodal points of the foreground mesh $\boldsymbol{x}_i$.
Using the (bi-) linearity of $a$ and $L$, we can re-write \eqref{eq:defn-k-f} as
\begin{equation}
K_{ij} = \sum_{k,\ell}\exMat_{ki}a(\phi_\ell,\phi_k)\exMat_{\ell j}\quad\text{and}\quad F_i = \sum_k \exMat_{ki}L(\phi_k)\text{ .}
\end{equation}
In matrix form, this reads as
\begin{equation}
    \mathbf{K} = \mathbf{\exMat}^T\mathbf{A}\mathbf{\exMat}\quad\text{and}\quad\mathbf{F} = \mathbf{\exMat}^T\mathbf{B}\text{ ,}\label{eq:mtam}
\end{equation}
where
\begin{equation}\label{eq:defn-a-b}
    A_{ij} = a(\phi_j,\phi_i)\quad,\quad B_i = L(\phi_i)\text{ ,}
\end{equation}
We will refer to $\mathbf{\exMat}$ as the ``extraction matrix''.  With some additional straightforward bookkeeping, this method extends to linear PDE systems with multiple scalar solution fields, as presented for extraction-based IGA in \cite{Kamensky2019}.  It can also be applied to each step of Newton iteration to solve a nonlinear problem, where $a$ is the Gateaux derivative of a residual form $-L$.

Crucially, the matrix $\mathbf{A}$ and vector $\mathbf{B}$ can be assembled using a standard boundary-fitted FE code in a loop over elements in the foreground mesh.  Thus, given a foreground mesh and an extraction matrix $\mathbf{\exMat}$, one can directly reuse an existing FE code for interpolation-based immersed analysis without making any modifications to the assembly procedures.  The post processing features of existing FE codes can also be reused by leveraging a foreground FE representation of the solution:
\begin{equation}
    u^h = \sum_{j=1}^\nu c_j\phi_j\text{ ,}
\end{equation}
where the foreground FE coefficients $\{c_i\}$ are obtained from the solution $\mathbf{d}$ of \eqref{eq:linear-system} by
\begin{equation}\label{eq:cmd}
    \mathbf{c} = \mathbf{\exMat}\mathbf{d}\text{ .}
\end{equation}

\subsection{Example implementation using FEniCS}
\label{subsec:fenics_impl}

For the numerical results in this paper, we apply the method outlined in Section \ref{subsec:exMat} to the existing FE automation toolchain FEniCS \cite{Logg2012}.  FEniCS allows arbitrary linear and bilinear forms to be specified over a variety of classical FE function spaces using the Python-based Unified Form Language (UFL) \cite{Alnaes2014}.  UFL forms are compiled into efficient element-level routines \cite{Kirby2006} which are used by the FE library DOLFIN \cite{Logg:2010} to assemble matrices and vectors.  In the notation of Section \ref{subsec:exMat}, UFL would be used to specify the forms $a$ and $L$, and DOLFIN would be used to assemble the matrix $\mathbf{A}$ and vector $\mathbf{B}$, applying the compiled element-level routines to elements $\{\omega_e\}_{i=1}^{\nu_\text{el}}$ of the foreground mesh.  DOLFIN's application--programmer interface (API) provides access to the linear algebra objects $\mathbf{A}$ and $\mathbf{B}$ as data structures in the parallel linear algebra library PETSc \cite{petsc-web-page,petsc-user-ref,petsc-efficient}, which is used to form $\mathbf{K}$ and $\mathbf{F}$ via \eqref{eq:mtam}, solve the linear system \eqref{eq:linear-system} for $\mathbf{d}$, and compute the foreground FE coefficients $\mathbf{c}$ via \eqref{eq:cmd}.  DOLFIN's output and visualization routines can then be applied to this foreground FE representation of the solution.

This non-invasive reuse of FEniCS provides great flexibility for implementing formulations of a variety of different PDE systems as we illustrate in the sequel.  
What remains is to generate the foreground- and background meshes alongside the extraction operators $\mathbf{M}$. The tools and methods used to do so are outlined in the following subsection.

\subsection{Mesh generation}
\label{subsec:mesh_generation}

This section describes the techniques used in the numerical examples of Sections \ref{sec:NumericalResults} and \ref{sec:Application} to generate both background-fitted and -unfitted foreground meshes.

\subsubsection{Background-fitted foreground meshes}
\label{subsec:bg_fitted}

\begin{figure}[t]
	\begin{center}
	\graphicspath{{}}
	\def\svgwidth{11.0cm}
\begingroup%
  \makeatletter%
  \providecommand\color[2][]{%
    \errmessage{(Inkscape) Color is used for the text in Inkscape, but the package 'color.sty' is not loaded}%
    \renewcommand\color[2][]{}%
  }%
  \providecommand\transparent[1]{%
    \errmessage{(Inkscape) Transparency is used (non-zero) for the text in Inkscape, but the package 'transparent.sty' is not loaded}%
    \renewcommand\transparent[1]{}%
  }%
  \providecommand\rotatebox[2]{#2}%
  \newcommand*\fsize{\dimexpr\f@size pt\relax}%
  \newcommand*\lineheight[1]{\fontsize{\fsize}{#1\fsize}\selectfont}%
  \ifx\svgwidth\undefined%
    \setlength{\unitlength}{113.2545315bp}%
    \ifx\svgscale\undefined%
      \relax%
    \else%
      \setlength{\unitlength}{\unitlength * \real{\svgscale}}%
    \fi%
  \else%
    \setlength{\unitlength}{\svgwidth}%
  \fi%
  \global\let\svgwidth\undefined%
  \global\let\svgscale\undefined%
  \makeatother%
  \begin{picture}(1,0.33136661)%
    \lineheight{1}%
    \setlength\tabcolsep{0pt}%
    \put(0,0){\includegraphics[width=\unitlength,page=1]{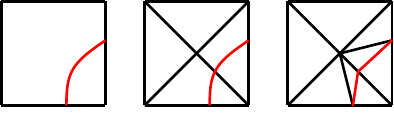}}%
    \put(0.11933345,0.01238927){\color[rgb]{0,0,0}\makebox(0,0)[lt]{\lineheight{1.25}\smash{\begin{tabular}[t]{l}(a)\end{tabular}}}}%
    \put(0.48229055,0.01065768){\color[rgb]{0,0,0}\makebox(0,0)[lt]{\lineheight{1.25}\smash{\begin{tabular}[t]{l}(b)\end{tabular}}}}%
    \put(0.85052758,0.01385886){\color[rgb]{0,0,0}\makebox(0,0)[lt]{\lineheight{1.25}\smash{\begin{tabular}[t]{l}(c)\end{tabular}}}}%
    \put(0.1392091,0.19653041){\color[rgb]{1,0,0}\makebox(0,0)[lt]{\lineheight{1.25}\smash{\begin{tabular}[t]{l}$\partial\Omega$\end{tabular}}}}%
  \end{picture}%
\endgroup%

	\caption{Generation of background fitted foreground meshes. (a) Single background element with description of the boundary shown by a red curve; (b) regular subdivision of the background element into triangles; (c) subdivision of triangles based on intersections of the interface with edges; interface approximated by red lines.} 
	\graphicspath{{Figures/}}
	\label{fig:mesh_decomp}
    \end{center}
\end{figure}
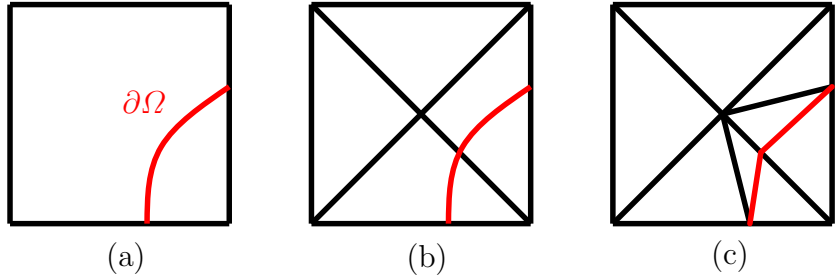

To generate background-fitted foreground meshes and their extraction matrices $\mathbf{M}$, the XFEM Tool Kit of the research code presented in \cite{Noel2022} was used. Within this toolkit, the foreground elements are created by repeated subdivision of a structured background grid. 
First, each background cell is subdivided into a set of triangular or tetrahedral elements as is shown in Figure \ref{fig:mesh_decomp} (b). Then, all edges of the elements are checked for intersections with the geometry. Using a set of subdivision templates, elements with intersected edges are then again subdivided into straight-edged triangular and tetrahedral elements such that new element vertices are created at the intersection locations. The new element facets (i.e. edges in 2D and faces in 3D) are a first-order geometric approximation of the domain boundary as shown in Figure \ref{fig:mesh_decomp} (c). 

The XFEM tool kit constructs the extraction matrices $\mathbf{M}$ such that basis functions are enriched following a Heaviside enrichment strategy \cite{Noel2022}. Enrichment is needed when the support of a basis contains two or more geometrically disconnected domains. Note that for the numerical examples studied in this paper enrichment is not required to ensure optimal convergence.

\vspace{0.5cm}

\begin{figure}[h!]
	\centering
	 \begin{subfigure}[b]{0.35\linewidth}
    \includegraphics[width=0.95\linewidth]{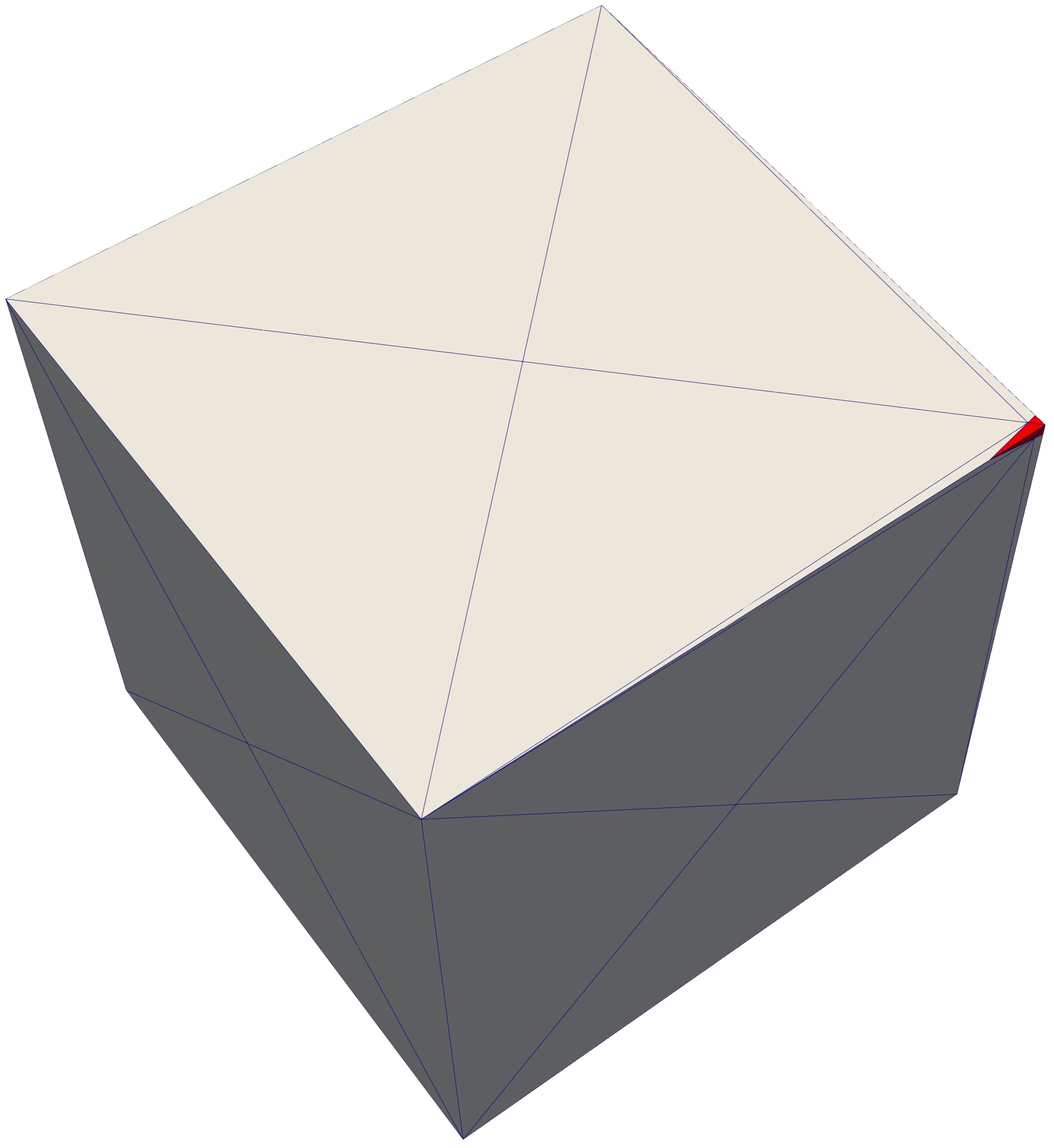} \caption{}\label{fig:sliver_element_cube}
  \end{subfigure}
  	 \begin{subfigure}[b]{0.35\linewidth}
    \includegraphics[width=0.95\linewidth]{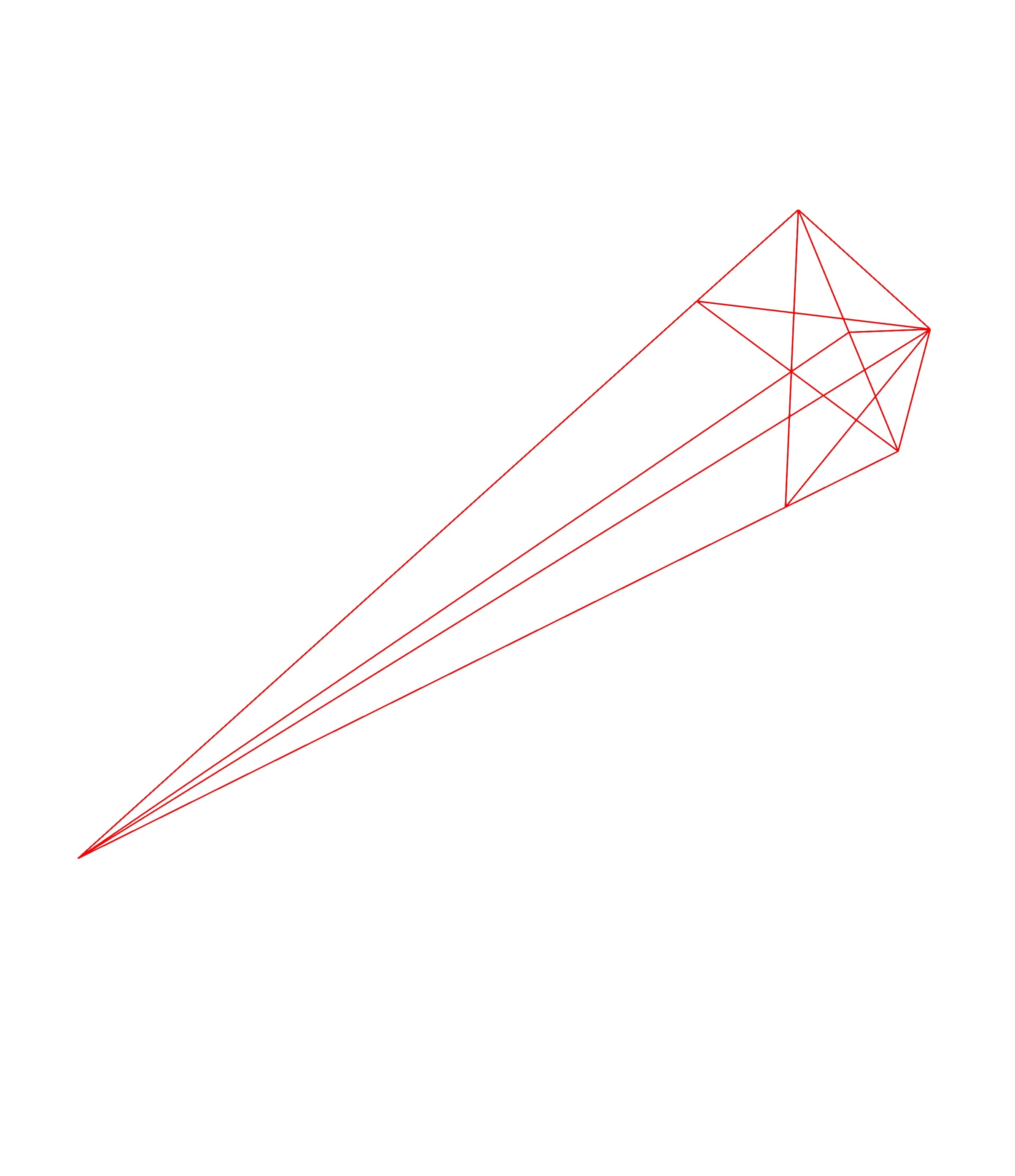}\caption{}\label{fig:sliver_element_detail}
  \end{subfigure}
    \caption{ (a) A tessellation of a hexahedral background element from the 3D mesh shown in Figure \ref{fig:poisson}, highlighting the sliver element in the rightmost corner. (b) Only the red section consisting of small tetrahedral elements intersects the PDE domain $\Omega$.} 
    \label{fig:sliver_element}
\end{figure}

While the procedure outlined enables the robust background-fitted triangulation of the material domain, it may lead to elements of arbitrarily small size and arbitrarily large aspect ratio.
As an example, the 2D foreground meshes used in Subsection \ref{subsec:bending_tab} contain elements with aspect ratios as high as $\sim 300$ and elements with volumes as small as $\sim 10^{-4}$ of the biggest elements. 
The poor mesh quality is exacerbated in 3D. The 3D foreground meshes used for the example in Subsection \ref{subsec:poisson} contain elements that differ by a factor of $\sim 10^{27}$ in volume, and aspect ratios as high as $\sim 10^{29}$. An example for the poor elements resulting from this mesh generation process in 3D is shown in Figure \ref{fig:sliver_element}.

\subsubsection{Background-unfitted foreground meshes}
\label{subsec:bg_unfitted}

Mesh cutting tools, like the XFEM tool kit discussed in Section \ref{subsec:bg_fitted}, generate background-fitted meshes and do not require an additional mesh generator to create foreground meshes. However, background-unfitted analysis allows substantially greater flexibility in the generation of foreground meshes; essentially any existing mesh generator can be used, with the only necessary information from the background mesh being an approximate element size.  We apply two different approaches within the examples shown in this paper, with the particular details largely motivated by convenience.  The background-unfitted foreground meshes used in Section \ref{subsec:unfitted} were generated using FEniCS's built in structured mesh functions, as were the corresponding background meshes. For these examples, the background function spaces were $C^0$ Lagrange FE spaces, and the extraction operators were generated directly in FEniCS, using its {\tt PETScDMCollection} functionality. 
The background-unfitted foreground meshes used in Sections \ref{subsec:pinned_shell} and \ref{subsec:bending_tab} use an isogeometric background basis. These meshed were created with FEniCS's mshr module, which generates unstructured simplicial meshes from geometric descriptions. The extraction operators were generated using the FEniCS-based isogeometric analysis library tIGAr \cite{Kamensky2019} by simply overriding methods that normally generate structured meshes for exact Lagrange extraction.

\section{Numerical results}
\label{sec:NumericalResults}

We now explore the properties of interpolation-based immersed analysis numerically using a variety of benchmark problems.  

\subsection{Comparing quadrature-based and interpolation-based methods} 
\label{subsec:poisson}

The first question we seek to answer is:  How do interpolation-based immersed methods with background-fitted foreground meshes compare with the more thoroughly-studied class of quadrature-based immersed methods?  (We defer the comparison of background-fitted and -unfitted foreground meshes to Section \ref{subsec:unfitted}.)  In particular, we use the Poisson equation as a model problem and consider a background function space of tensor-product B-splines of degree $k$ (at maximal continuity) and a foreground Lagrange FE space of equal polynomial degree $\kappa = k$ on a simplicial mesh.  The background function space contains monomials of degree greater than $k$ which cannot be represented exactly by the foreground function space, so interpolation-based methods differ from quadrature-based ones.  

To carry out this test, we construct an instance of the Poisson problem using the method of manufactured solutions, where we choose the exact solutions
\begin{equation}
    u(\boldsymbol{x}) = \sin(\pi(x_1^2 + x_2^2))\cos(\pi(x_1 - x_2))
\end{equation}
and 
\begin{equation}
    u(\boldsymbol{x}) = \sin(\pi(x_1^2 + x_2^2 + x_3^2))\cos(\pi(x_1 + x_2 + x_3))
\end{equation}
in 2D and 3D respectively, and define the source term as $f := -\Delta u$.  For the 2D problem, $\Omega$ is taken to be a unit square rotated by $45^\circ$ about the origin (i.e., $\Omega = \{\boldsymbol{x}\in\mathbb{R}^2~:~\Vert\boldsymbol{x}\Vert_{\ell^1}<1/2\}$), and the manufactured solution $u$ is applied as Dirichlet data on the boundary of this rotated square. $\Omega$ is then immersed into a structured background mesh of an axis-aligned bi-unit square on which the background B-spline spaces are constructed.  A suite of progressively-refined pairs of background and background-fitted foreground meshes are used with background element size $h = 2^{-(R+1)}$, where $R\in\{0,\cdots,6\}$ denotes the refinement level and ``$h$'' is the square root of background element area. 

For the 3D case, we take $\Omega$ to be an analogous rotated unit cube (rotated by $45^\circ$ about the $x_3$ axis and then by $45^\circ$ about the $x_2$ axis) and define a sequence of structured hexahedral background meshes with element size $h = 2^{-(R+1)}$, where $R\in\{0,\cdots,4\}$ is the refinement level and $h$ is the cube root of the element volume.  The foreground meshes in 3D are background-fitted tetrahedral meshes.  

\begin{figure}[!t]
	\centering
	 \begin{subfigure}[b]{0.45\linewidth}
    \includegraphics[width=\linewidth]{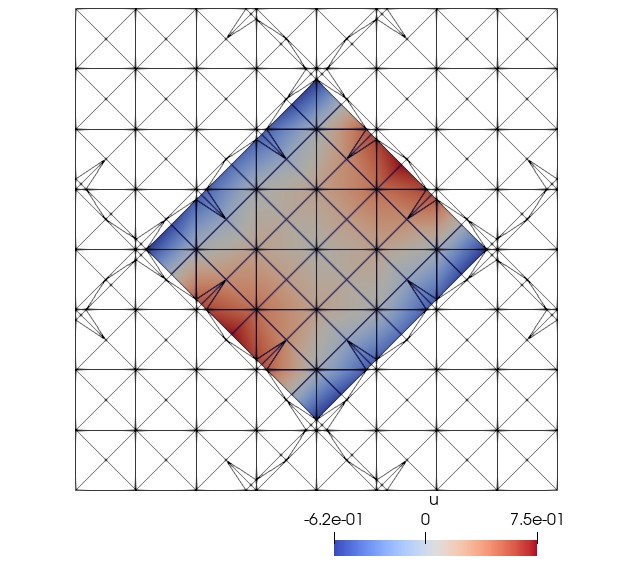} \caption{\label{fig:poisson2d}}
  \end{subfigure}
  	 \begin{subfigure}[b]{0.5\linewidth}
    \includegraphics[width=\linewidth]{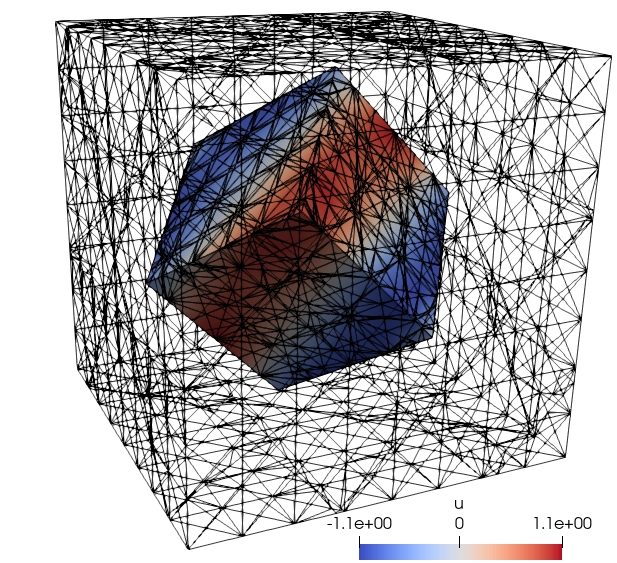}\caption{\label{fig:poisson3d}}
  \end{subfigure}
    \caption{Numerical solution to the Poisson problem using interpolation-based immersed analysis, visualized on the foreground mesh of refinement level $R=2$. Tests were performed on both the 2D square (a) and 3D rotated cube (b). 
    } 
    \label{fig:poisson}
\end{figure}

Dirichlet boundary conditions are enforced weakly using the non-symmetric Nitsche formulation \eqref{eq:poisson-disc-generic}, with $C_\text{pen} = 0$.  A representative foreground mesh and numerical solution are shown in Figure \ref{fig:poisson}.

\begin{figure}[t]
	\centering
    \includegraphics[width=\linewidth]{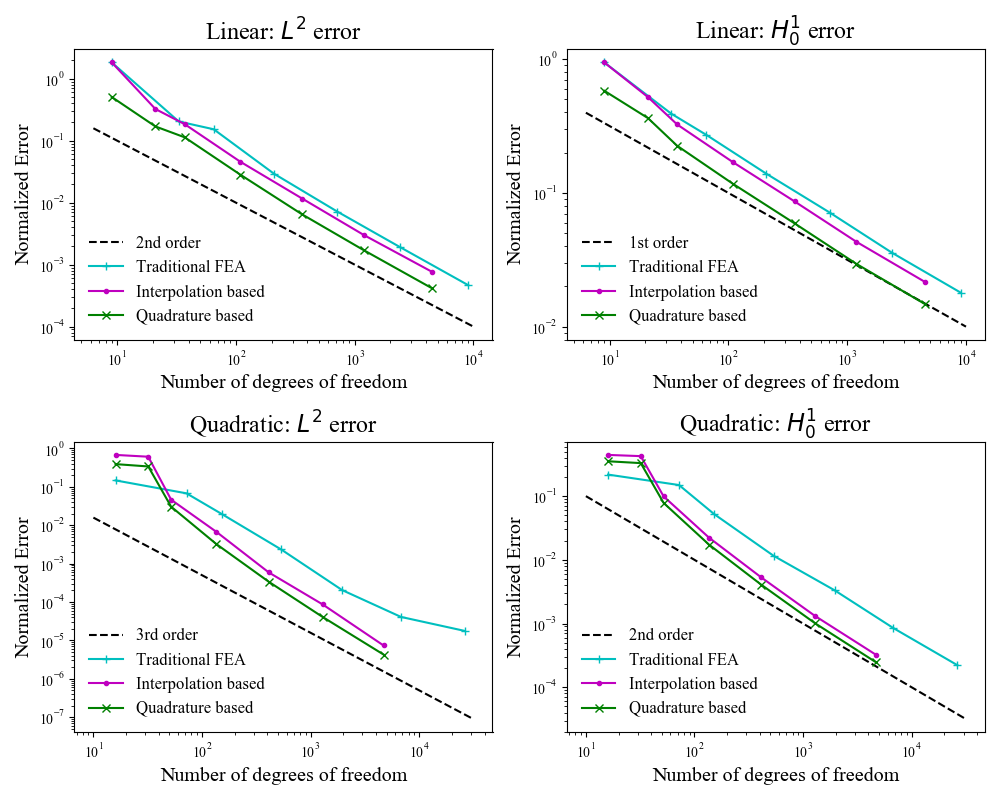}
    \caption{\label{fig:poisson_plots2D} Convergence data for Poisson's problem in 2D, comparing  the results of traditional FEA on the foreground mesh with the Lagrange FE space, the interpolation based immersed method developed in this paper, and the quadrature based immersed method described in Section \ref{sec:model-quadrature-based}.}
\end{figure}

\begin{figure}[t]
\centering
    \includegraphics[width=\linewidth]{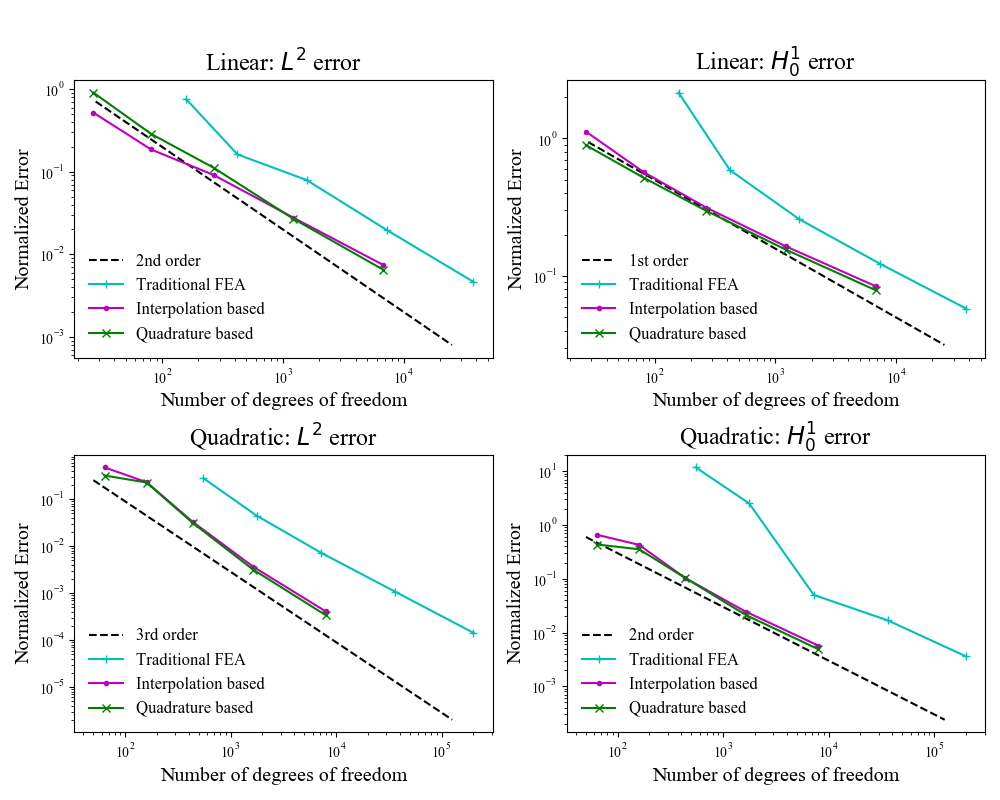}
    \caption{\label{fig:poisson_plots3D} Convergence data for Poisson's problem in 3D, comparing  the results of traditional FEA on the foreground mesh with the Lagrange FE space, the interpolation based immersed method developed in this paper, and the quadrature based immersed method described Section \ref{sec:model-quadrature-based}. }
\end{figure}

Figures \ref{fig:poisson_plots2D} and \ref{fig:poisson_plots3D} compare the convergence of interpolation-based and quadrature-based immersed calculations.\footnote{The interpolation-based calculations use the FEniCS-based implementation described in Section \ref{sec:Implementation}, whereas the quadrature-based calculations use a different research code with the requisite functionality which is outlined in \cite{Noel2022} }  Both types of immersed method are seen to converge optimally for $k=1$ and $k=2$ in both $H^{1}_{0}$ and $L^2$ norms.  We also include convergence results of standard FE analysis using the Lagrange FE space defined on the foreground mesh.  This also converges optimally with respect to $h$, but is less accurate per degree of freedom due to excess foreground mesh refinement in cut background elements and less-efficient approximation of smooth functions than the maximally-smooth B-splines of the background function space.  

\subsection{Non-conforming background spaces} 
\label{subsec:biharmonic}

The numerical analysis sketched in Section \ref{sec:Approximation} made the assumption of a conforming discrete space, i.e., that the discrete space $\mathcal{V}^h$ of interpolated background functions was a subset of the infinite dimensional space used to define the continuous variational problem (viz., $H^1$ in the case of the Poisson equation).  However, one may wish to relax that requirement in practice to eliminate topological constraints on the foreground mesh and/or interpolate smooth background functions inexactly with less-regular foreground FE functions. 

As a model problem to illustrate non-conforming interpolation-based immersed analysis, we consider the biharmonic problem:  Find $u$ such that
\begin{equation}
    \Delta^2 u = f\text{ ,}
\end{equation}
with boundary conditions
\begin{align}
    u = \sigma &\text{ on } \partial\Omega\text{ ,} \\
    \grad u \cdot \boldsymbol{n} = \grad \sigma \cdot \boldsymbol{n} &\text{ on } \partial\Omega\text{ ,}
\end{align}
where $f:\Omega\to\mathbb{R}$ is a given source term and $\sigma:\Omega\to\mathbb{R}$ is an auxiliary function used to define Dirichlet boundary data for $u$ and $\nabla u\cdot\boldsymbol{n}$ on $\partial\Omega$.\footnote{One could instead define separate functions to prescribe boundary data for $u$ and its normal derivative, but an equivalent $\sigma$ must exist for the problem to be well-posed and the present perspective is more convenient for applying the method of manufactured solutions.}  Unlike the nonsymmetric Nitsches method used in the Poisson problem in \eqref{eq:poisson-disc-generic}, we discretize this problem using a symmetric Nitsche-like method:  Find $u^h\in\mathcal{V}^h$ such that, $\forall v^h\in\mathcal{V}^h$, 
\begin{align} 
    \nonumber \int_{\Omega}\Delta u^h\Delta v^h\,d\Omega + \int_{\partial\Omega}\grad \Delta u^h \cdot \boldsymbol{n}v^h - \Delta u^h \grad v^h \cdot  \boldsymbol{n}\,d\Gamma & \\
    \nonumber  + \int_{\partial\Omega}(\grad \Delta v^h) \cdot \boldsymbol{n}(u^h - \sigma)  - \Delta v(\grad u^h \cdot  \boldsymbol{n} - \grad \sigma  \cdot  \boldsymbol{n})\,d\Gamma  & \\
    + \int_{\partial\Omega}\dfrac{\alpha}{h^3} (u^h - \sigma)v^h + \dfrac{\beta}{h} (\grad u^h \cdot  \boldsymbol{n} - \grad \sigma  \cdot  \boldsymbol{n})\grad v^h\cdot  \boldsymbol{n}\,d\Gamma  & = \int_\Omega fv^h\,d\Omega\text{ ,}\label{eq:biharmonic-nitsche}
\end{align}
where $\alpha > 0$ and $\beta > 0$ are independent of $h$ and $eta$, but must be sufficiently large to ensure stability of the formulation.\footnote{For the method to be usable with piecewise-quadratic function spaces, $\alpha$ must be strictly positive, even if a non-symmetric variant is used, as $\nabla\Delta v^h \equiv 0$ and only the penalty term is left to enforce $u^h\approx\sigma$.} In the computations of this paper, we use $\alpha=\beta=5$.  In the case of immersed methods, the integrals $\int_\Omega(\cdots)\,d\Omega$ and $\int_{\partial\Omega}(\cdots)\,d\Gamma$ should be understood as being split up over elements of the foreground mesh, as in \eqref{eq:split-integrals}.  This distinction is important when non-conforming spaces are used as some spatial derivatives may have distributional contributions at interior mesh facets that are ignored, so \eqref{eq:biharmonic-nitsche} technically involves some abuse of notation.

As with the Poisson problem in Section \ref{subsec:poisson}, we use the method of manufactured solutions, choosing exact solutions
\begin{equation}
    u(\boldsymbol{x}) = \cos(0.05\pi x_1+0.1)\cos(0.05\pi x_2 +0.1))
\end{equation}
and 
\begin{equation}
    u(\boldsymbol{x}) =\cos(\pi x_1+0.5)\cos(\pi x_2 +0.5)\cos(\pi x_3 + 0.5)
\end{equation}
in 2D and 3D, respectively, and setting the problem data to $f := \Delta^2 u$ and $\sigma := u$.

\begin{figure}[t]
	\centering
	\includegraphics[width= \textwidth]{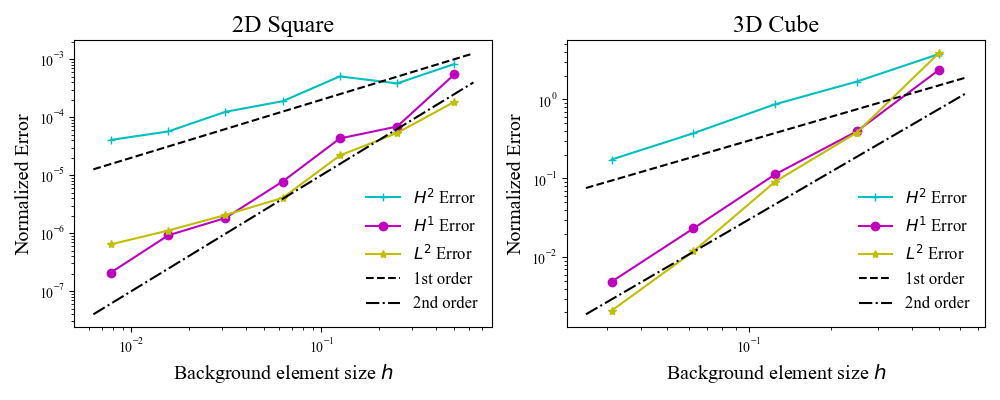}\
	\caption{\label{fig:biharmonic}Convergence data for the biharmonic problem.}
\end{figure}

We reuse the rotated square and cube domains and mesh sequences from Section \ref{subsec:poisson}.  We select the {\em background} function spaces to be $H^2$-conforming.  In particular, we use quadratic B-splines of maximal continuity.  However, the {\em foreground} function spaces are taken to be quadratic Lagrange FE spaces, which, on simplicial meshes, do not represent all monomials of the background B-spline spaces.  Further, interpolations of background functions on this Lagrange FE space are only $C^0$, containing jumps in their first derivatives (so they can be at most $H^{3/2-\epsilon}$-conforming, for $\epsilon>0$ arbitrarily-small).  Interpolation-based immersed analysis will therefore be non-conforming in this case.  2D and 3D convergence results are shown in Figure \ref{fig:biharmonic}.  We get optimal rates of 1 and 2 for the $H^2$ and $H^1$ norms respectively, where ``$H^2$ norm'' is understood in a broken sense (i.e., the square root of a sum over squared foreground-element $H^2$ norms) in this non-conforming setting.  We get only a rate of 2 for the $L^2$ norm, but this is expected, since an Aubin--Nitsche-type duality argument to obtain an optimal $L^2$ error bound breaks down for quadratic B-spline discretizations of the biharmonic problem \cite{StrFix73}.  

\begin{remark}\label{rem:volume-filter}
The 3D foreground meshes used in this section include some elements where $\eta \ll h$, such as the example shown in Figure \ref{fig:sliver_element}.  The harmful effects of small background element--domain intersections on linear system conditioning have been studied previously in the context of quadrature-based immersed methods, along with a variety of possible remedies, but we leave the adaptation of such techniques to the interpolation-based setting to future work, as discussed further in Section \ref{sec:conclusion}.  However, small foreground elements, where $\eta\to 0$, present a distinct challenge unique to interpolation-based methods, due to the blow-up of foreground basis function derivatives .\footnote{This blow-up does not manifest in quadrature-based methods since they directly evaluate background basis function derivatives.}  The resulting large numbers in the matrix $\mathbf{A}$ defined by \eqref{eq:defn-a-b} should, in exact arithmetic, be canceled by small entries in the extraction matrix $\mathbf{M}$ when forming $\mathbf{K}$ by \eqref{eq:mtam}.  However, this cancellation is not exact due to finite precision.  For the fourth-order biharmonic problem in 3D, we find that this can prevent convergence even with direct solvers.  In the 3D biharmonic examples from the present work, we circumvent the issue by removing near-zero volume elements from the assembly of $\mathbf{A}$, skipping integration over foreground elements (and their associated boundary facets) if their volume falls below $0.001\%$ of the largest foreground element.  The resulting geometric perturbation of the domain has no observable effect on solutions, as illustrated by the convergence results of Figure \ref{fig:biharmonic}.  One can imagine pathological domain--mesh intersections where such a strategy would lead to significant errors, but these do not arise in our tests, and the solution of fourth-order problems on topologically-3D domains is itself relatively uncommon in applications.
\end{remark}

\subsection{Geometric approximations of domain geometry } 
\label{subsec:linearElasticity}

The convergence tests of Sections \ref{subsec:poisson} and \ref{subsec:biharmonic} deliberately used polyhedral PDE domains so that a foreground mesh of affine simplex elements could fit them exactly, allowing us to isolate errors due to interpolation of the background basis functions.  We now consider a problem with a curved boundary where there is some geometry error between $\Omega$ and the foreground mesh.  

In particular, we choose the classic linear elasticity benchmark of an infinite 2D plate with a circular hole subjected to biaxial tension.  We consider the plane strain variant of the problem.  An exact solution for the stress field is readily available in elasticity textbooks, e.g., \cite[Section 7.7.5]{Lubarda2020}.  We truncate the problem to a finite square domain whose side length is four times the diameter of the hole and apply the traction from the exact solution to the infinite-domain problem on the boundary.  This finite-domain problem is further truncated by restricting it to the upper-right quadrant of the $x_1$--$x_2$ plane and applying symmetry boundary conditions on the axes, which eliminates rigid-body modes from the displacement solution of the original pure Neumann problem.  We formulate linear elasticity as a problem for displacement, viz., Navier's equations, which also provide an opportunity to demonstrate interpolation-based immersed analysis for a system of PDEs:  Find a displacement field $\boldsymbol{u}:\Omega\to\mathbb{R}^2$ such that
\begin{equation}
    -\nabla\cdot\boldsymbol{\sigma}(\boldsymbol{u}) = \boldsymbol{0}
\end{equation}
subject to boundary conditions
\begin{align}
    \boldsymbol{u}\cdot\boldsymbol{n} = 0\quad&\text{on}~\Gamma_\text{sym}\text{ ,}\\
    (\boldsymbol{I}-\boldsymbol{n}\otimes\boldsymbol{n})(\boldsymbol{\sigma}\cdot\boldsymbol{n}) = \boldsymbol{0}\quad&\text{on}~\Gamma_\text{sym}\text{ ,}\\
    \boldsymbol{\sigma}\cdot\boldsymbol{n} = \boldsymbol{t}\quad&\text{on}~\Gamma_t\text{ ,}
\end{align}
where $\Gamma_\text{sym}\subset\partial\Omega$ is the union of the symmetry planes where sliding boundary conditions are applied, $\Gamma_t = \partial\Omega\setminus\Gamma_\text{sym}$ is the Neumann boundary on which traction is prescribed, $\boldsymbol{t}$ is the traction on $\Gamma_t$ taken from the exact solution, and the Cauchy stress $\boldsymbol{\sigma}$ is 
\begin{equation}
    \boldsymbol{\sigma}(\boldsymbol{u}) = 2\mu\boldsymbol{\varepsilon}(\boldsymbol{u}) + \lambda\operatorname{tr}(\boldsymbol{\varepsilon}(\boldsymbol{u}))\boldsymbol{I}
\end{equation}
in terms of the symmetric gradient operator $\boldsymbol{\varepsilon}(\cdot)$ (which gives strain when applied to the displacement) and Lam\'{e} parameters $\mu$ and $\lambda$, which we set from Young's modulus $E=200\times 10^9$ and Poisson's ratio $\nu = 0.3$ by standard formulas.

We discretize Navier's equations with slip boundary conditions using a symmetric Nitsche-like formulation analogous to \eqref{eq:poisson-disc-generic}:  Find $\boldsymbol{u}^h\in\mathcal{V}^h$ such that, $\forall\boldsymbol{v}^h\in\mathcal{V}^h$,
\begin{align}
    \nonumber &\int_\Omega \boldsymbol{\sigma}(\boldsymbol{u}^h):\nabla\boldsymbol{v}^h\,d\Omega\\
    \nonumber &\quad - \int_{\Gamma_\text{sym}}(\boldsymbol{u}^h\cdot\boldsymbol{n})\boldsymbol{n}\cdot\boldsymbol{\sigma}(\boldsymbol{v}^h)\cdot\boldsymbol{n}\,d\Gamma - \int_{\Gamma_\text{sym}} (\boldsymbol{v}^h\cdot\boldsymbol{n})\boldsymbol{n}\cdot\boldsymbol{\sigma}(\boldsymbol{u}^h)\cdot\boldsymbol{n}\,d\Gamma \\
    &\quad + \int_{\Gamma_\text{sym}}\frac{\beta \mu}{h}(\boldsymbol{u}^h\cdot\boldsymbol{n})\boldsymbol{n}\cdot\boldsymbol{v}^h\,d\Gamma = \int_{\Gamma_t}\boldsymbol{t}\cdot\boldsymbol{v}^h\,d\Gamma\text{ ,}
    \end{align}
where $\beta > 0$ is a penalty parameter associated with Nitsche's method (in this work we set $\beta = 10$) and integrals should be understood in the sense of \eqref{eq:split-integrals} (which is technically an abuse of notation in the present case because it ignores the geometry error in the foreground mesh).  Similarly, the outward-facing unit normal $\boldsymbol{n}$ depends on the geometry of the foreground mesh, but this dependence has been suppressed in the notation.

\begin{figure}
	\centering
	\begin{subfigure}[b]{0.32\textwidth}
	\centering
	 \includegraphics[width=\textwidth]{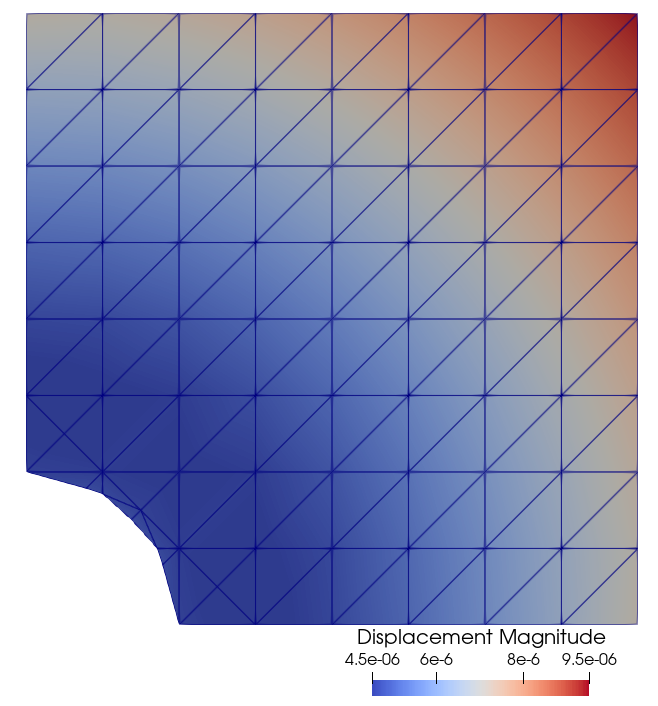} 
	  \caption{No foreground refinement }
	\end{subfigure}
	\hfill
	\begin{subfigure}[b]{0.32\textwidth}
	\centering
	 \includegraphics[width=\textwidth]{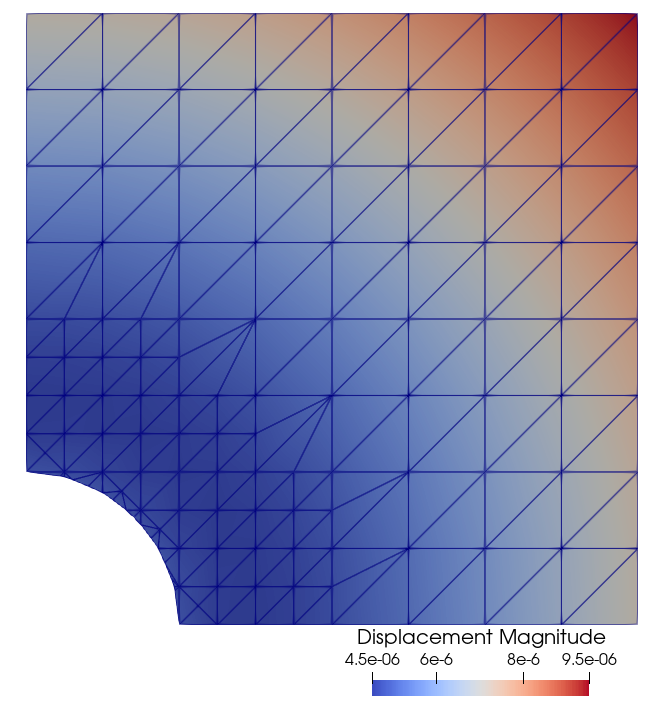} 
	  \caption{One level of foreground refinement }
	\end{subfigure}
	\hfill
	\begin{subfigure}[b]{0.32\textwidth}
	\centering
	 \includegraphics[width=\textwidth]{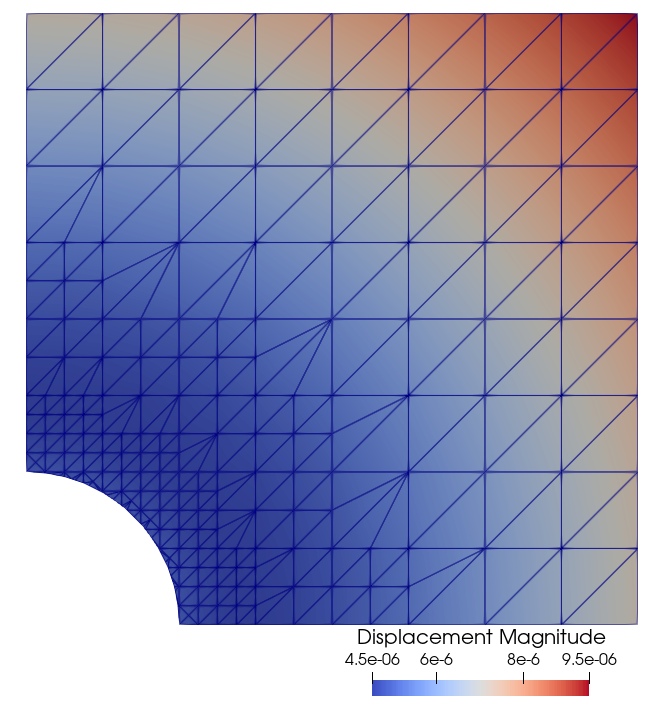} 
	  \caption{Two levels of foreground refinement }
	\end{subfigure}
	\caption{\label{fig:holeInPlate} Plots showing the magnitude of displacement. To resolve curved surfaces, the foreground mesh is locally refined, while the same uniform background mesh is used.}
\end{figure}

To study the behavior of interpolation-based immersed methods, we construct a sequence of uniform, structured background meshes of quadrilateral elements, covering the truncated plate and hole.  We consider several corresponding sequences of foreground meshes, with different levels of local refinement around the curved boundary of the hole.  The three different foreground mesh refinement levels are shown for an $8\times 8$ background mesh in Figure \ref{fig:holeInPlate}, along with representative displacement solutions.  In all cases, the foreground mesh elements are affine triangles.  For the background function spaces, we select B-spline spaces of degrees $k=1$ and $k=2$, with maximal continuity for each, and, for the foreground function spaces, we use Lagrange FE spaces of corresponding polynomial degrees $k = \kappa$.  

\begin{figure}[h!]
	\centering
	\includegraphics[width= 0.5\textwidth]{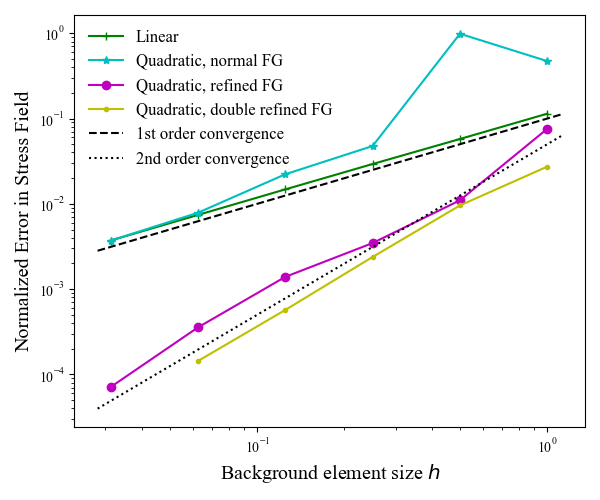}
	\caption{\label{fig:holeInPlateGraph} Frobenius norm of stress error tensor plotted against the background element size shown with both linear basis functions and quadratic basis functions. In the case of quadratic basis functions, the foreground mesh was refined as shown in Figure \ref{fig:holeInPlate}.}
\end{figure}

Figure \ref{fig:holeInPlateGraph} shows the convergence of the $L^2$ norm of error in the Cauchy stress for different polynomial degrees and levels of foreground refinement.  We see that, even without foreground refinement the case of $k=1$ attains the optimal convergence rate.  For $k=2$, geometry error dominates in the case without foreground refinement and limits the rate of convergence.  However, even with just one level of foreground refinement near the curved boundary, we see an optimal effective convergence rate.  For affine foreground elements and a fixed level of local foreground refinement, we still expect the geometry error to dominate in the ultimate asymptotic limit (as occurs also with the even lower-order stair-step approximations of curved boundaries in finite cell adaptive quadrature), but this phenomenon is unlikely to be encountered in practical problems where the model error would dominate well before this asymptotic regime.  We hypothesize that the use of high-order isoparametric elements instead of local foreground refinement would provide true optimal convergence, but high-order foreground mesh generation is currently outside the capabilities of the software framework used for this paper.

\subsection{Nonlinear and unsteady problems} 
\label{subsec:TGVortex}
We now consider the unsteady incompressible Navier--Stokes equations as an example of a more complicated nonlinear PDE system.  This problem also illustrates the use of interpolation-based immersed finite element analysis in conjunction with stabilized finite element formulations that use mesh-dependent variational forms, originally designed with the standard body-fitted FE setting in mind.  

In particular, we consider the PDE system:  Find velocity and pressure fields $(\boldsymbol{u},p):\Omega\times (0,T)\to\mathbb{R}^d\times\mathbb{R}$ such that
\begin{align}
    \nonumber \rho \left(\partial_t\boldsymbol{u} + \boldsymbol{u}\cdot \grad \boldsymbol{u}\right)  &= - \grad p + \mu \Delta \boldsymbol{u} \\
    \grad \cdot \boldsymbol{u } &= 0 
\end{align}
subject to the boundary condition
\begin{equation}
    \boldsymbol{u} = \boldsymbol{g}\quad\text{on}~\partial\Omega\text{ ,}
\end{equation}
where $d$ is the spatial dimension of $\Omega$, $(0,T)$ is the time interval on which the problem is posed, $\rho$ is mass density, $\mu$ is dynamic viscosity, and $\boldsymbol{g}$ is given Dirichlet boundary data.  We discretize this in space using a stabilized variational multiscale (VMS) formulation \cite{Bazilevs2007a} with non-symmetric Nitsche-like weak enforcement of the Dirichlet boundary condition \cite{Bazilevs07c}, leading to the semi-discrete problem:  Find $(\boldsymbol{u}^h,p^h)\in\pmb{\mathcal{V}}^h\times\mathcal{Q}^h$ such that, $\forall (\boldsymbol{v}^h,q^h)\in\pmb{\mathcal{V}}^h\times\mathcal{Q}^h$,
\begin{align}
    \nonumber\int_\Omega \left( \rho \left( \partial_t\boldsymbol{u}^h + \boldsymbol{u}^h \cdot \grad \boldsymbol{u}^h\right) \cdot \boldsymbol{v}^h  +\boldsymbol{\sigma}(\boldsymbol{u}^h,p^h):\nabla\boldsymbol{v}^h  +  \grad \cdot \boldsymbol{u }^h q^h  \right)\,d\Omega &\\ 
    \nonumber  + \int_{\partial \Omega} \left(- (\boldsymbol{\sigma}(\boldsymbol{u}^h,p^h)\boldsymbol{n})\cdot\boldsymbol{v}^h  +  (\boldsymbol{\sigma}(\boldsymbol{v}^h,q^h) \cdot \boldsymbol{n}) \cdot (\boldsymbol{u}^h - \boldsymbol{g})\right)\,d\Gamma\\
    \nonumber - \int_{\partial\Omega}\left(\rho \min\{\boldsymbol{u}^h\cdot\boldsymbol{n},0\}(\boldsymbol{u}^h - \boldsymbol{g}) \cdot \boldsymbol{v}^h \right) \,d\Gamma&\\
    \nonumber + \int_\Omega \tau_{\text{M}}\left(  \left( \boldsymbol{u}^h  \cdot \grad \boldsymbol{v}^h + \dfrac{1}{\rho}\grad q^h \right)\cdot \boldsymbol{r}_{\text{M}}   - \boldsymbol{v}^h \cdot \boldsymbol{r}_{\text{M}} \cdot \grad \boldsymbol{u}^h\right)\,d\Omega&\\
    \nonumber - \int_\Omega\dfrac{\tau_{\text{M}}^2}{\rho} (\grad \boldsymbol{v}^h):( \boldsymbol{r}_{\text{M}}\otimes \boldsymbol{r}_{\text{M}})\,d\Omega& \\
    + \int_\Omega \tau_{\text{C}} \left(r_{\text{C}} \grad \cdot \boldsymbol{v}^h \right)\,d\Omega   &= 0 \label{eq:vms}
\end{align}
where 
\begin{equation}
    \boldsymbol{\sigma}(\boldsymbol{u},p) = 2\mu\boldsymbol{\varepsilon}(\boldsymbol{u}) - p\boldsymbol{I}
\end{equation}
gives the functional form of the Cauchy stress in terms of the velocity and pressure fields, 
\begin{equation}
    \boldsymbol{r}_\text{M} = \rho (\partial_t\boldsymbol{u}^h + \boldsymbol{u}^h  \cdot \grad \boldsymbol{u}^h) -  \nabla\cdot\boldsymbol{\sigma}(\boldsymbol{u}^h,p^h)
\end{equation} 
is the residual of the strong momentum balance equation, 
\begin{equation}
r_{\text{C}} = \rho \grad \cdot \boldsymbol{u}^h    
\end{equation}
is the residual of the strong continuity equation, and $\tau_{\text{M}}$ and $\tau_{C}$ are the momentum and continuity stabilization parameters, respectively.  A number of definitions for these stabilization parameters have been studied in the literature.  
We follow \cite{Bazilevs08a} in defining 
\begin{equation}
    \tau_{\text{M}} = \left(\boldsymbol{u}^h \cdot \boldsymbol{G} \cdot \boldsymbol{u}^h + C_{\text{I}} \nu^{2} \boldsymbol{G}:\boldsymbol{G} + \frac{C_{t}}{\Delta t^2} \right)^{- 1/2},
\end{equation}
and 
\begin{equation}
    \tau_{\text{C}} = \dfrac{1}{\tau_{\text{M}} \text{tr}(\mathbf{G})},
\end{equation}
with $\nu = \mu/\rho$ the kinematic viscosity, $\Delta t$ the time step size of the temporal discretization, $C_\text{I} > 0$ and $C_t > 0$ are dimensionless constants (set to $C_\text{I} = 60$ and $C_t=4$ in this work), and $\boldsymbol{G}$ is an anisotropic mesh size tensor, which simplifies to
\begin{equation}
    \boldsymbol{G} = 4 h^{-2}\boldsymbol{I}
\end{equation}
in the case of a uniform isotropic mesh with element size $h$.  In the computations of this paper, we take this element size to be from the uniform background mesh.  We discretize the system of ordinary differential equations emanating from the semi-discrete problem \eqref{eq:vms} using the implicit midpoint rule in time, with a time step $\Delta t$ proportional to the background element size $h$.  

As a representative instance of the Navier--Stokes problem, we consider the 2D Taylor--Green vortex.  The statement and spatially-periodic exact solution of this problem are well-known.  In numerical studies, it is often restricted to the square $(0,\pi)^2$ with symmetry boundary conditions (e.g., \cite[Section 9.10.1]{Evans2011}).  However, in the present work, we take $\Omega$ to be the rotated unit square used as the domain of the Poisson problem in Section \ref{subsec:poisson}, apply the exact solution as the boundary data $\boldsymbol{g}$, and use the same sequence of background and foreground mesh pairs defined in Section \ref{subsec:poisson}.  The stabilized formulation \eqref{eq:vms} permits the stable use of equal-order discretizations, i.e.,
\begin{equation}
    \pmb{\mathcal{V}}^h = (\mathcal{V}^h)^d\quad\text{and}\quad\mathcal{Q}^h = \mathcal{V}^h
\end{equation}
for a single scalar function space $\mathcal{V}^h$.  Our numerical tests in this section use the same linear and quadratic choices for $\mathcal{V}^h$ as used for the scalar Poisson problem in Section \ref{subsec:poisson}.  

\begin{figure}[t]
	\centering
	\includegraphics[width= 1\textwidth]{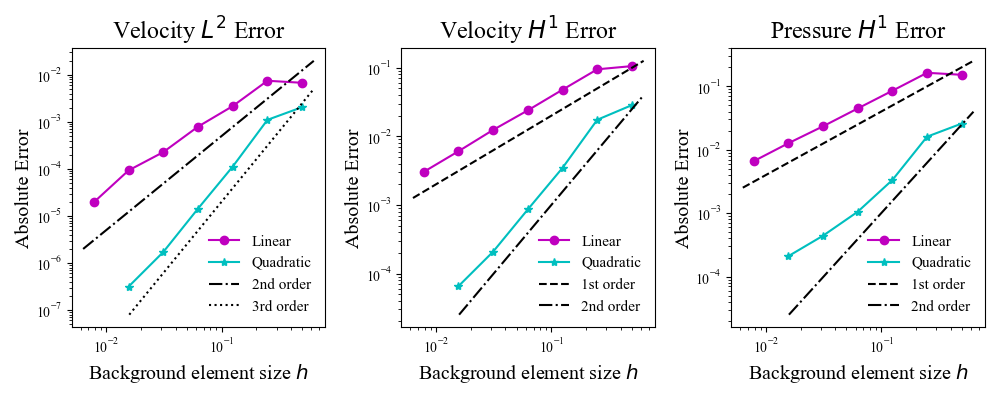}
	\caption{\label{fig:TG} Convergence data for the velocity and pressure fields of the Taylor--Green vortex problem using both linear and quadratic discretizations. Kinematic viscosity was set to $1/100$ yielding a Reynolds number of 100,  and error is evaluated at time $T=1$.}
\end{figure}

The convergence of velocity and pressure as $h\sim\Delta t \to 0$ (resulting in a CFL number of approximately $1$) is shown in Figure \ref{fig:TG}. We see effectively optimal rates of velocity convergence with respect to the polynomial degree of the spatial discretizations in the $L_2$ and $H^1$ norms.  In principle, using $h\sim\Delta t$ and the implicit midpoint rule should limit convergence to second-order in the $L_2$ norm,\footnote{The assumption of $h\sim\Delta t$ is tacit in the stabilized formulation used; refining $\Delta t$ faster than $h$ would weaken stabilization in the asymptotic regime.  See \cite{Hsu2010} for additional discussion.} but the quadratic spatial discretizations appear to be third-order in velocity $L^2$ error within the range of resolutions tested.  Pressure $H^1$ error convergence is also nearly optimal, but with a drop in convergence rate on the finest meshes when using quadratic spatial discretizations.  This sub-optimal $H^1$ pressure convergence is consistent with the analysis of related stabilized methods for Stokes and Oseen flow, e.g., \cite{Hughes1986,Franca1992}, where error is bounded in a norm on the velocity--pressure product space in which the pressure gradient is scaled by element size, weakening its contribution.

\subsection{Background-unfitted foreground meshes}
\label{subsec:unfitted}

The analysis of Section \ref{sec:Approximation} and preceding numerical examples in Sections \ref{subsec:poisson}--\ref{subsec:TGVortex} all used background-fitted foreground meshes.  Background-fitted foreground meshes are a necessity for quadrature-based immersed methods, since the quadrature rules defined on foreground elements require smoothness of the integrand for accuracy.  However, interpolations of background basis functions on a foreground mesh are, by construction, smooth within each foreground element and can be accurately integrated with foreground element-by-element quadrature rules.  Of course, the solution space for interpolation-based immersed methods consists of interpolations of background basis functions that are non-smooth within each foreground element, so the analysis of Section \ref{sec:Approximation} is not applicable.  However, the numerical results of the current section suggest that optimal accuracy is maintained with background-unfitted foreground meshes.  

\begin{figure}[t]
	\centering
	\includegraphics[width= 0.5\textwidth]{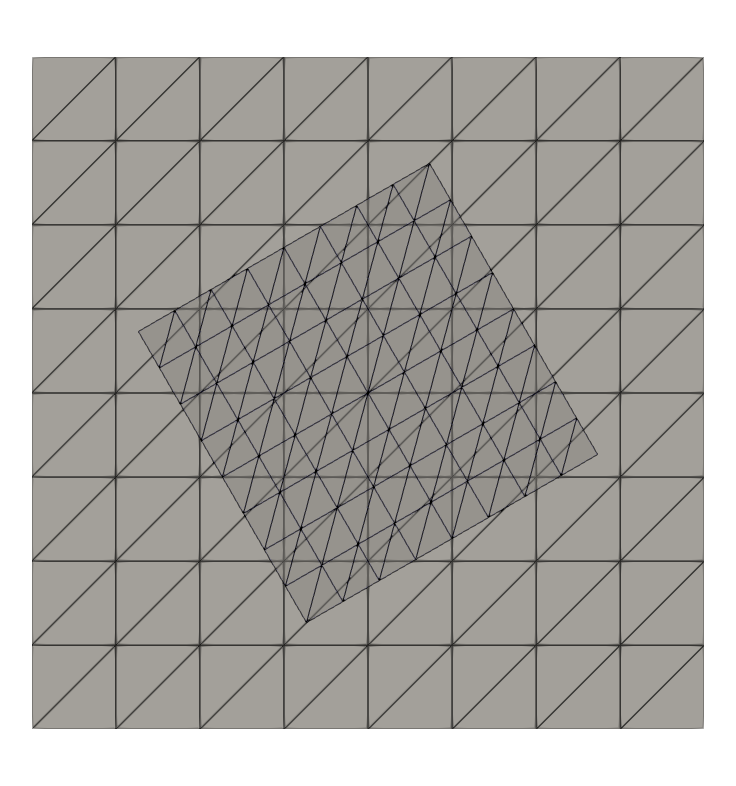}
	\caption{\label{fig:unfittedMeshes} A representative example of a background mesh and a background-unfitted foreground mesh.}
\end{figure}

In the preceding examples, the background function space was taken to be a B-spline space of maximal continuity.  However, as a stress test of background-unfitted interpolation-based immersion, we instead test it with $C^0$ Lagrange background function spaces defined on simplicial meshes.  As test problems, we consider the Poisson equation with a manufactured solution (cf.\ Section \ref{subsec:poisson}) and the Taylor--Green vortex (cf.\ Section \ref{subsec:TGVortex}).  For both of these test problems, the domain is taken to be a rotated square, and the sequence of pairs of background and background-unfitted foreground meshes are constructed similarly to the representative pair shown in Figure \ref{fig:unfittedMeshes} but at varying levels of refinement.  In all tests, the foreground element size $\eta$ is kept proportional to the background element size $h$, and the polynomial degrees $k$ and $\kappa$ of the background and foreground are equal.  For the Taylor--Green vortex problem, the same time implicit midpoint rule with time step $\Delta t$ proportional to the background element size $h$ was used as in Section \ref{subsec:TGVortex}. 

\begin{figure}[h!]
	\centering
	\includegraphics[width= 0.8\textwidth]{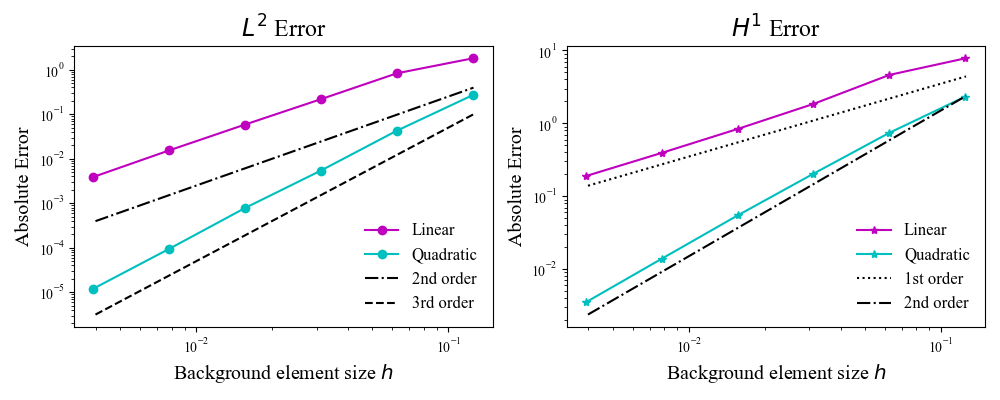}
	\caption{\label{fig:unfittedPoisson}  Convergence data for the Poisson problem, generated using extraction onto  unfitted background meshes (Figure \ref{fig:unfittedMeshes}).}
\end{figure}

\begin{figure}[h!]
	\centering
	\includegraphics[width= 0.8\textwidth]{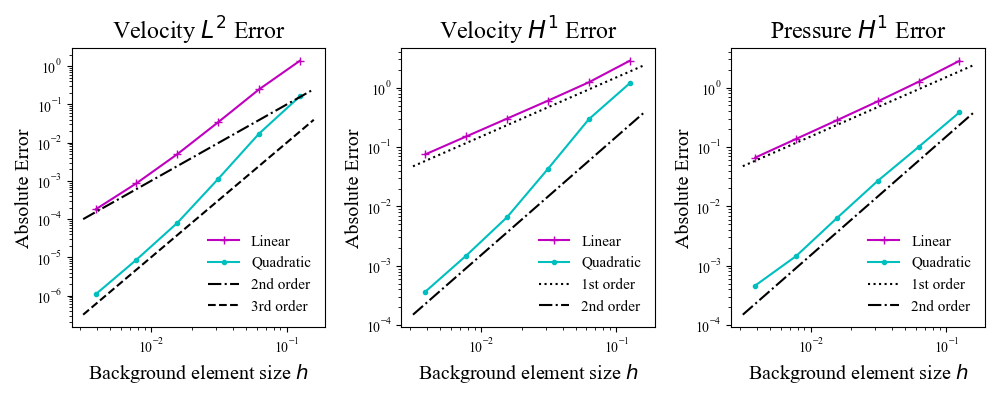}
	\caption{\label{fig:unfittedTG} Convergence data for the Taylor--Green vortex problem  generated using extraction onto unfitted meshes (Figure \ref{fig:unfittedMeshes}), with Reynolds number 100, after elapsed time $T=1$. }
\end{figure}

The convergence results for both the Poisson and Navier--Stokes examples are shown in Figures \ref{fig:unfittedPoisson} and \ref{fig:unfittedTG} respectively.  We see optimal convergence rates of $L^2$ and $H^1$ error (of velocity in the case of Navier--Stokes) for linear and quadratic discretizations.

\section{Application to IMGA of trimmed shells}
\label{sec:Application}
A major success of IGA has been its application to thin shell analysis, originating in the work of Kiendl and collaborators \cite{Kiendl2009,Kiendl2011,Kiendl2015}.  Smooth discrete function spaces from IGA make it possible to apply the Bubnov--Galerkin method directly to Kirchhoff--Love theory where the displacement satisfies a fourth-order PDE system.  This provides much better per-degree-of-freedom accuracy than classical FE discretizations of shells, sometimes resolving important qualitative solution features with orders of magnitude less computation \cite{Morganti:2015fy}.  However, the original isogeometric Kirchhoff--Love shell formulation of \cite{Kiendl2009} is restricted to boundary-fitted spline representations of the midsurface geometry.  With suitably-flexible spline technologies, such as T-splines \cite{Seder03a}, it is practical to design such ``analysis-suitable'' geometries for industrial applications as demonstrated for thin shell analysis of composite wind turbine blades \cite{BazHsu12b}, prosthetic heart valve leaflets \cite{Hsu2015}, and automobile parts \cite{Casquero2020}.  However, the vast majority of industrial CAD software defines surface geometries in an immersed way by cutting them out of unfitted parametric coordinate charts along arbitrary trim curves \cite{Marussig2018}.  Thus, immersogeometric formulations are needed to directly analyze CAD models of shell structures.  Quadrature-based immersed methods for thin shells have previously been used in a number of studies \cite{Guo2018,Coradello2020,Coradello2021}.  In this section of the present work, we apply interpolation-based IMGA to Kirchhoff--Love analysis of shells whose geometries are defined by subsets of a B-spline patch's parameter space.  

The specific variational formulation that we use in this section is that of a geometrically-nonlinear Kirchhoff--Love shell of uniform thickness $h_\text{th}$ whose material behavior follows a homogeneous isotropic St.\,Venant--Kirchhoff model defined by a Young's modulus $E$ and Poisson ratio $\nu$.  A complete mathematical statement of this problem is given by \cite{Kiendl2011}, and its translation into FEniCS UFL is available online in the library ShNAPr \cite{shnapr-code} developed to support the work in \cite{Kamensky2021}.  We therefore omit the extensive kinematic definitions for brevity as they are neither novel nor essential to following the present study.  We mention briefly, though, that we implement pinned boundary conditions on the edges of trimmed shells using a penalty formulation based on \cite{herrema2019penalty}.  This consists of adding
\begin{equation}
    +\frac{1}{2}\int_{\mathcal{L}_\text{pin}}\frac{\alpha Eh_{\text{th}}}{h}\left\vert\boldsymbol{u}^h - \boldsymbol{g}\right\vert^2\,d\mathcal{L}\label{eq:penalty-energy}
\end{equation}
to the shell structure's elastic energy functional where $\mathcal{L}_\text{pin}$ is the portion of the trim curve (in the reference configuration) on which pinned boundary conditions are to be applied, $\alpha > 0$ is a dimensionless constant, $h$ is the background-mesh element size, $\boldsymbol{u}^h$ is the unknown discrete displacement field of the shell structure's midsurface, and $\boldsymbol{g}$ is the prescribed midsurface displacement on the pinned boundaries.  For quadratic spline discretizations, this ``naive'' penalty method for pinned boundaries is essentially as accurate as a full Nitsche-type formulation (cf.\ \cite{Benzaken2021}) because the associated consistency term involves third-order derivatives.  The use of \eqref{eq:penalty-energy} in conjunction with an analogous rotational penalty has been demonstrated to be accurate for coupling non-matching shell structures in practical problem settings \cite{herrema2019penalty,Johnson2020,Liu2021,Zhao2022}.

\subsection{Verification with boundary-fitted IGA} 
\label{subsec:pinned_shell}

\begin{figure}
     \centering
     \begin{subfigure}[b]{0.45\textwidth}
         \centering
         \includegraphics[width=\textwidth]{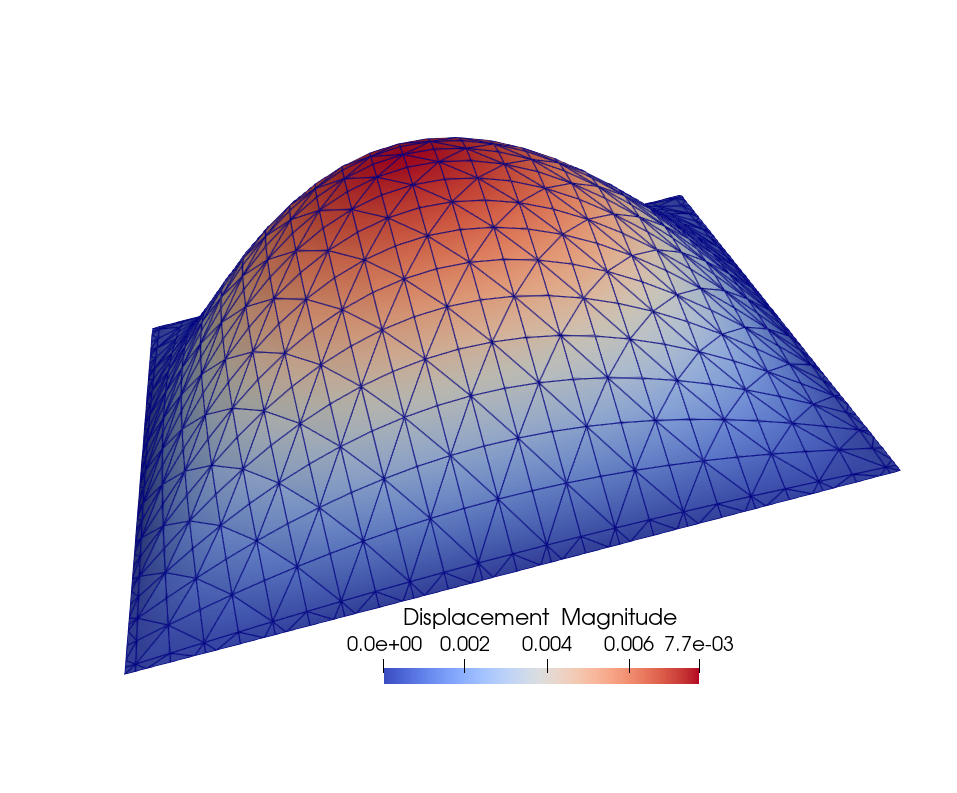}
         \caption{}\label{fig:fitted_plate}
     \end{subfigure}
     \hfill
     \begin{subfigure}[b]{0.45\textwidth}
         \centering
         \includegraphics[width=\textwidth]{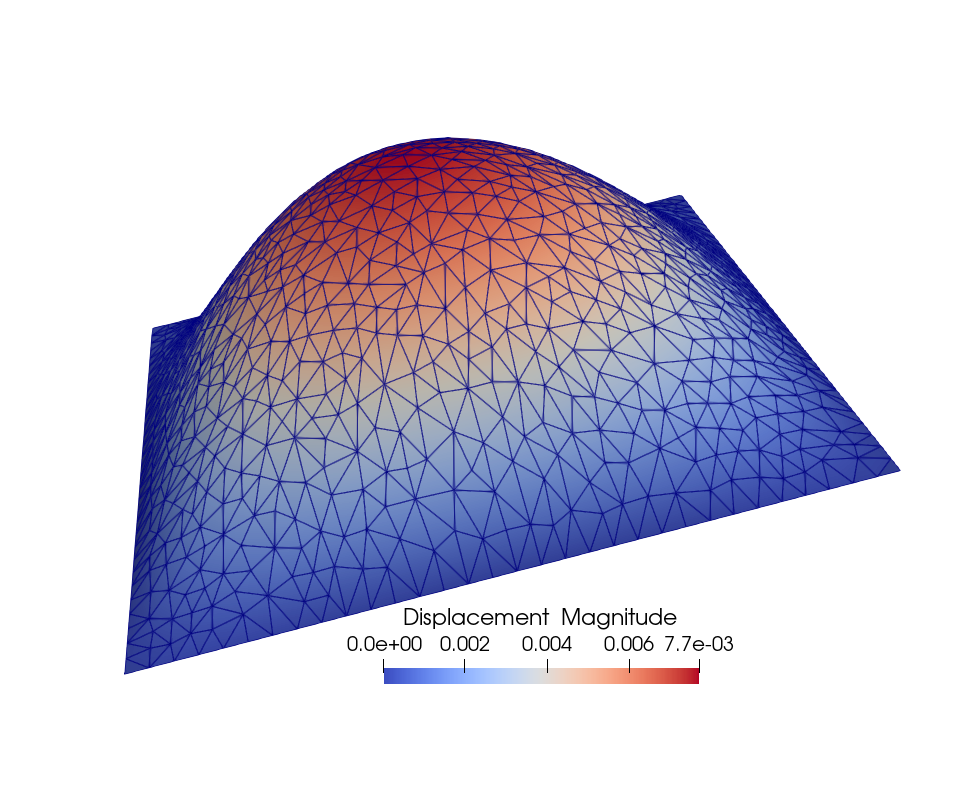}
         \caption{}\label{fig:unfitted_plate}
     \end{subfigure}
    \hfill
     \begin{subfigure}[b]{0.5\textwidth}
         \centering
         \includegraphics[width=\textwidth]{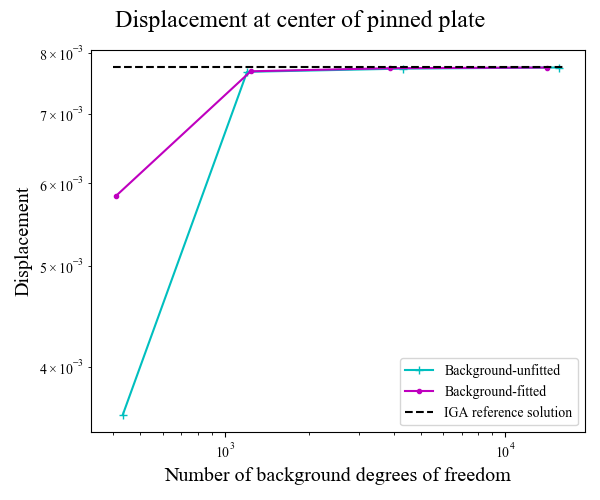}
         \caption{}\label{fig:pinned_plot}

     \end{subfigure}
\caption{\label{fig:pinned_plate} Examples of the background-fitted (a) and background-unfitted (b) discretizations used to verify the IMGA formulation. (c) Displacement of the center of the plate is plotted against the number of degrees of freedom in each mesh.. }
\end{figure}

We first test interpolation-based IMGA of thin shells by trimming out a rectangular domain which can easily be modeled with an untrimmed B-spline to provide a reference solution for verification purposes.  In particular, for each Cartesian component of the shell midsurface displacement, we reuse the $k=2$ sequence of scalar B-spline background and Lagrange-interpolated foreground spaces from the 2D Poisson and biharmonic tests of Sections \ref{subsec:poisson} and \ref{subsec:biharmonic}, where the parametric-space domain is a rotated unit square trimmed out of an axis-aligned bi-unit square, shown in Figure \ref{fig:pinned_plate}. Additionally, we use a new background-unfitted $k=2$ sequence of scalar B-spline background and Lagrange-interpolated foreground spaces, constructed as described in Section \ref{subsec:bg_unfitted}. In this test, we take the mapping from the B-spline parameter space to the physical-space reference configuration of the shell's midsurface to be given by 
\begin{equation}
    \boldsymbol{X}(\boldsymbol{\xi}) = \xi_1\boldsymbol{e}_1 + \xi_2\boldsymbol{e}_2\text{ ,}
\end{equation}
where $\boldsymbol{\xi} = (\xi_1,\xi_2)$ is a point in the B-spline parameter space, $\boldsymbol{X}$ is the corresponding 3D point on the reference configuration midsurface, and $\{\mathbf{e}_1,\mathbf{e}_2,\mathbf{e}_3\}$ is the standard orthonormal basis for a Cartesian coordinate system of 3D physical space into which the reference configuration is embedded.  The construction of an analogous sequence of untrimmed reference B-spline discretizations is straightforward due to the square geometry.

We apply homogeneous (i.e., $\boldsymbol{g} = \boldsymbol{0}$) pinned boundary conditions to the entire boundary of the unit square domain, using penalty coefficient $\alpha = 10^5$, and assume material parameters of $E = 4.8\times 10^5$ and $\nu = 0.38$, and thickness $h_\text{th} = 0.1$.  We then apply a uniform force density per unit reference area of $90\boldsymbol{e}_3$, leading to an ``inflated'' equilibrium configuration, of which a representative numerical approximation using interpolation-based IMGA is shown in Figures \ref{fig:fitted_plate} and \ref{fig:unfitted_plate}.  We also extract the out-of-plane displacement component at the center of the plate using IMGA at several levels of refinement and plot it as a function of number of degrees of freedom in Figure \ref{fig:pinned_plot}.   This displacement converges to the reference value obtained from a highly-refined computation using standard untrimmed IGA.  (This is a geometrically-nonlinear problem, for which there is no analytical Kirchhoff--Love solution available.)
 
\subsection{Background-fitted and -unfitted IMGA for a complicated trimmed geometry}
\label{subsec:bending_tab}

\begin{figure}
     \centering
     \begin{subfigure}[b]{0.45\textwidth}
         \centering
         \includegraphics[width=\textwidth]{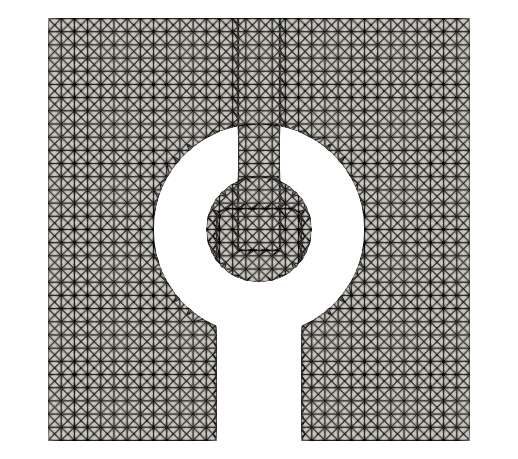}
         \caption{Background-fitted mesh}
         \label{fig:bent_par_f}
     \end{subfigure}
     \hfill
     \begin{subfigure}[b]{0.45\textwidth}
         \centering
         \includegraphics[width=\textwidth]{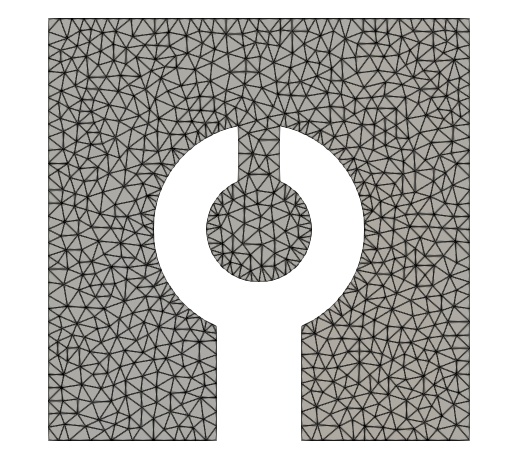}
         \caption{Background-unfitted mesh}
         \label{fig:bent_par_u}
     \end{subfigure}
     \hfill
     \begin{subfigure}[b]{0.45\textwidth}
         \centering
         \includegraphics[width=\textwidth]{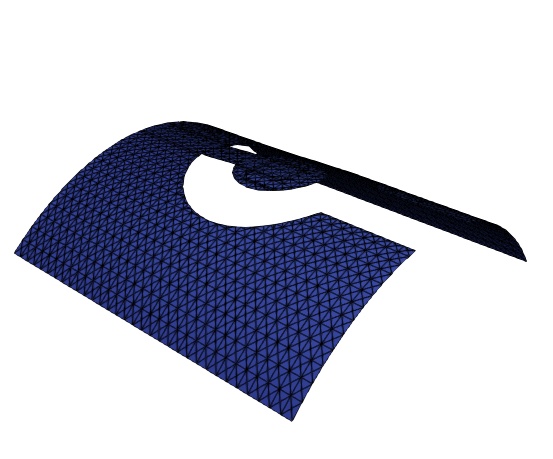}
         \caption{Background-fitted reference configuration}
         \label{fig:bent_ref_f}
     \end{subfigure}
     \begin{subfigure}[b]{0.45\textwidth}
         \centering
         \includegraphics[width=\textwidth]{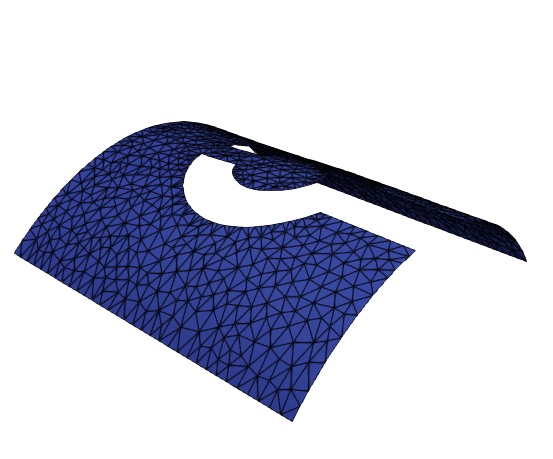}
         \caption{Background-unfitted reference configuration}
         \label{fig:bent_ref_u}
     \end{subfigure}
     \hfill
     \begin{subfigure}[b]{0.45\textwidth}
         \centering
         \includegraphics[width=\textwidth]{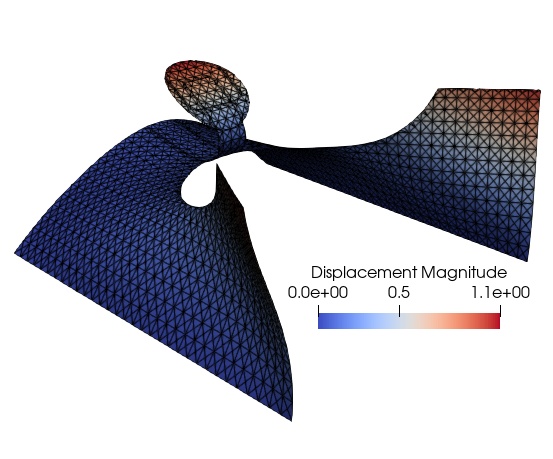}
         \caption{Background-fitted deformed configuration}
         \label{fig:bent_final_f}
     \end{subfigure}
     \hfill
     \begin{subfigure}[b]{0.45\textwidth}
         \centering
         \includegraphics[width=\textwidth]{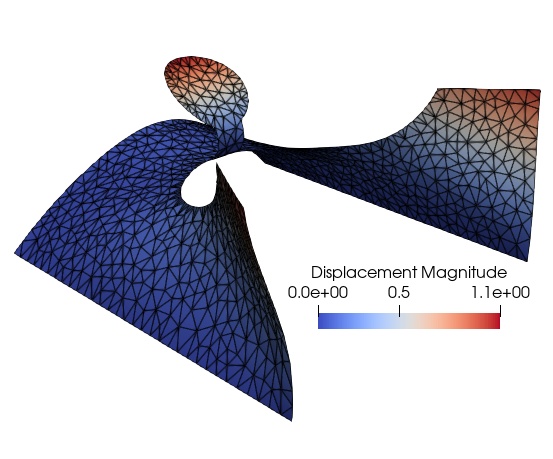}
         \caption{Background-unfitted deformed configuration}
         \label{fig:bent_final_u}
     \end{subfigure}
        \caption{Background-fitted and -unfitted meshes (in B-spline parameter space) and reference and deformed configurations (in physical space) of a trimmed shell structure subjected to a follower load.}
        \label{fig:bent}
\end{figure}

As a final example, we consider a more complicated trimmed spline geometry, for which a boundary-fitted IGA model would require unstructured spline technologies that are not currently part of standard design workflows.  Background-fitted and -unfitted meshes of the trimmed geometry are shown in the spline parameter space in Figures \ref{fig:bent_par_f} and \ref{fig:bent_par_u}. The parametric space is then mapped to a curved physical-space reference configuration by
\begin{equation}
    \boldsymbol{X}(\boldsymbol{\xi}) = \xi_1\boldsymbol{e}_1 + \xi_2\boldsymbol{e}_2 + A(1-\xi_1^2)\boldsymbol{e}_3\text{ ,}
\end{equation}
as visualized in Figures \ref{fig:bent_ref_f} and \ref{fig:bent_ref_u} where $A =1/2$.  The physical problem that we solve on this geometry is a static Kirchhoff--Love shell analysis with $E = 3\times10^{4}$, $\nu =0.3$, and $h_\text{th} =0.03 $, pinned boundary conditions applied to the outer left and right sides of the parametric domain, and a pressure follower load with force density per unit reference midsurface area of $\boldsymbol{f} = 90\boldsymbol{n}$ where $\boldsymbol{n}$ is the unit normal vector to the deformed midsurface.  Representative numerical approximations of the resulting equilibrium configuration are shown in Figures \ref{fig:bent_final_f} and \ref{fig:bent_final_u}.

\begin{remark}
The static equilibrium configurations shown in Figures \ref{fig:bent_final_f} and \ref{fig:bent_final_u} were arrived at by computing the steady-state solution of a dynamic Kirchhoff--Love shell problem with strong mass damping and numerical dissipation from under-resolved backward-Euler time stepping.  This pseudo-time continuation method was necessary to avoid bifurcations that arise on some meshes when attempting to apply static load stepping.  
\end{remark}

\begin{figure}[h!]
	\centering
	\includegraphics[width=0.5\textwidth]{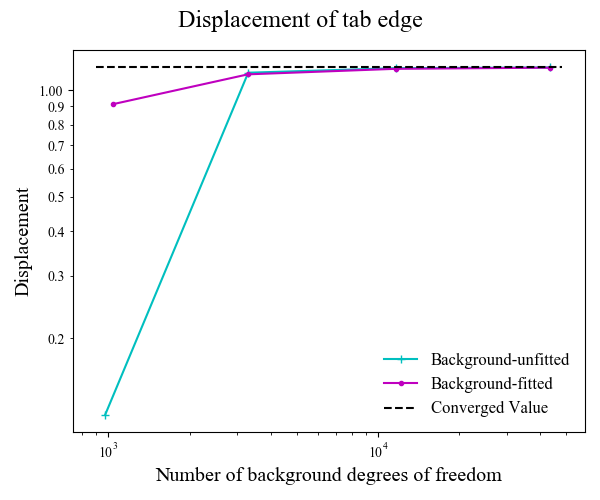}
	\caption{\label{fig:bent_plot} Convergence of background-fitted and -unfitted analyses of a trimmed shell subject to a follower load.}
\end{figure}

To study convergence with respect to mesh refinement quantitatively, we extract the displacement magnitude at the parametric point $\boldsymbol{\xi} = (0,-1/4)$ (corresponding to the bottom of the inner circle) as a quantity of interest.  This displacement is shown as a function of number of degrees of freedom for both the background-fitted and -unfitted discretizations in Figure \ref{fig:bent_plot}.   We can see that both types of discretizations converge to a consistent value as the number of degrees of freedom increases, where we estimate the converged value of $0.0078$ based on a highly-refined background-unfitted discretization.

\section{Conclusion}\label{sec:conclusion}

This paper introduces a new paradigm for defining and implementing high-order-accurate immersed-boundary methods based on interpolating basis functions defined on a boundary-unfitted background mesh into a function space defined on a foreground mesh that is fitted to the domain boundary but subject to fewer topological or geometrical constraints than a standard FE mesh.  Compared with existing variational immersed methods that use foreground meshes to define quadrature rules for functions on background meshes, interpolation-based immersed finite element analysis provides more flexibility with respect to the construction of the foreground mesh and more straightforward reuse of existing FE software.  The interpolation-based methodology is capable of exactly reproducing quadrature-based immersed methods under certain conditions, but it also enables efficient approximations that reduce computational cost with no significant loss of accuracy.  

One major challenge of quadrature-based immersed methods that remains in the interpolation-based paradigm is the problem that the conditioning of algebraic equation systems like \eqref{eq:linear-system} can become arbitrarily-bad in certain situations where a background mesh element's intersection with the PDE domain is very small \cite{DePrenter2017}.  In the present work, we have not applied any special treatment of these situations aside from the {\em ad hoc} volume based cell-filtering discussed in Remark \ref{rem:volume-filter} for the 3D biharmonic tests, although that was motivated by floating-point precision issues, not linear system conditioning.\footnote{The technique of Remark \ref{rem:volume-filter} may nonetheless incidentally improve system conditioning, by effectively removing basis functions whose supports have small intersections with the domain, similar to the approach analyzed by \cite{Elfverson2018}.}  In practice, we observe that iterative linear solvers remain effective for the examples considered in this paper, but we recognize that pathological conditioning like that identified in \cite{DePrenter2017} can still arise in interpolation-based immersed methods.  We hypothesize that techniques from the literature on quadrature-based immersed methods can be adapted {\em mutatis mutandis} to the interpolation-based setting with similar benefits, such as basis function agglomeration \cite{Kummer2016}, ghost penalty formulations \cite{Burman2010}, basis function removal \cite{Elfverson2018}, and specialized preconditioners \cite{Jomo2019}.  However, verifying this hypothesis is left to future work.  

While Section \ref{sec:Approximation} briefly sketches a path toward analyzing the approximation properties of the interpolated solution space, our verification efforts are primarily empirical in the present work, relying on numerical experiments with a variety of PDE systems.  These experiments indicate that interpolation-based immersed methods are optimally-accurate, even in conditions that are far outside the scope of our initial analysis.  It is therefore an interesting topic for future research to perform a more comprehensive {\em a priori} analysis of interpolation-based immersed methods, especially in the case of foreground meshes that are not fitted to elements of the background mesh.  

The non-invasive implementation of interpolation-based immersed methods discussed in Section \ref{sec:Implementation} is of great practical significance, as implementation cost has been a major barrier to widespread adoption and application of immersed methods proposed in the academic literature.  We have made the FEniCS-based prototype used for experiments in this paper available online \cite{interp-fea-code}.  However, we plan to develop a more comprehensive open-source tool chain for non-invasive implementation of interpolation-based immersed methods in the near future.  

\section*{Acknowledgements}
J. E. Fromm was supported by National Science Foundation award number 1650112.  R. Xiang and D. Kamensky were partially supported by National Science Foundation award number 2103939.  N. Wunsch, K. Maute, and J. A. Evans were supported by National Science Foundation award number 2104106.  H. Zhao was partially supported by National Aeronautics and Space Administration award number 80NSSC21M0070.

\bibliographystyle{unsrt}
\bibliography{exhumePaper1}

\begin{thebibliography}{10}

\bibitem{Hughes87}
T.~J.~R. Hughes.
\newblock {\em The Finite Element Method. Linear Static and Dynamic Finite
  Element Analysis}.
\newblock Prentice-Hall, Englewood Cliffs, New Jersey, 1987.

\bibitem{Hardwick2005}
M.~F. Hardwick, R.~L. Clay, P.~T. Boggs, E.~J. Walsh, A.~R. Larzelere, and
  A.~Altshuler.
\newblock {DART} system analysis.
\newblock Technical Report SAND2005-4647, Sandia National Laboratories, 2005.

\bibitem{Peskin72}
C.~S. Peskin.
\newblock Flow patterns around heart valves: A numerical method.
\newblock {\em Journal of Computational Physics}, 10(2):252--271, 1972.

\bibitem{Parvizian07}
J.~Parvizian, A.~D{\"u}ster, and E.~Rank.
\newblock Finite cell method.
\newblock {\em Computational Mechanics}, 41(1):121--133, 2007.

\bibitem{Schillinger2015}
D.~Schillinger and M.~Ruess.
\newblock The finite cell method: A review in the context of higher-order
  structural analysis of {CAD} and image-based geometric models.
\newblock {\em Archives of Computational Methods in Engineering},
  22(3):391--455, Jul 2015.

\bibitem{Burman2015}
E.~Burman, S.~Claus, P.~Hansbo, M.~G. Larson, and A.~Massing.
\newblock {CutFEM}: Discretizing geometry and partial differential equations.
\newblock {\em International Journal for Numerical Methods in Engineering},
  104(7):472--501, 2015.

\bibitem{Sulsky1995}
D.~Sulsky, S.-J. Zhou, and H.~L. Schreyer.
\newblock Application of a particle-in-cell method to solid mechanics.
\newblock {\em Computer physics communications}, 87(1-2):236--252, 1995.

\bibitem{Hughes05a}
T.~J.~R. Hughes, J.~A. Cottrell, and Y.~Bazilevs.
\newblock Isogeometric analysis: {CAD}, finite elements, {NURBS}, exact
  geometry, and mesh refinement.
\newblock {\em Computer Methods in Applied Mechanics and Engineering},
  194:4135--4195, 2005.

\bibitem{CoHuBa09}
J.~A. Cottrell, T.~J.~R. Hughes, and Y.~Bazilevs.
\newblock {\em Isogeometric Analysis: Toward Integration of CAD and FEA}.
\newblock Wiley, Chichester, 2009.

\bibitem{Marussig2018}
B.~Marussig and T.~J.~R. Hughes.
\newblock A review of trimming in isogeometric analysis: Challenges, data
  exchange and simulation aspects.
\newblock {\em Archives of Computational Methods in Engineering},
  25(4):1059--1127, Nov 2018.

\bibitem{Kamensky2015}
D.~Kamensky, M.-C. Hsu, D.~Schillinger, J.~A. Evans, A.~Aggarwal, Y.~Bazilevs,
  M.~S. Sacks, and T.~J.~R. Hughes.
\newblock An immersogeometric variational framework for fluid--structure
  interaction: Application to bioprosthetic heart valves.
\newblock {\em Computer Methods in Applied Mechanics and Engineering},
  284:1005--1053, 2015.

\bibitem{Elfverson2018}
D.~Elfverson, M.~G. Larson, and K.~Larsson.
\newblock {CutIGA} with basis function removal.
\newblock {\em Advanced Modeling and Simulation in Engineering Sciences},
  5(1):6, Mar 2018.

\bibitem{Johansson2019}
A.~Johansson, B.~Kehlet, M.~G. Larson, and A.~Logg.
\newblock Multimesh finite element methods: Solving {PDE}s on multiple
  intersecting meshes.
\newblock {\em Computer Methods in Applied Mechanics and Engineering},
  343:672--689, 2019.

\bibitem{Logg2012}
A.~Logg, K.-A. Mardal, G.~N. Wells, et~al.
\newblock {\em Automated Solution of Differential Equations by the Finite
  Element Method}.
\newblock Springer, 2012.

\bibitem{BSEH11}
M.~J. Borden, M.~A. Scott, J.~A. Evans, and T.~J.~R. Hughes.
\newblock Isogeometric finite element data structures based on {B}\'{e}zier
  extraction of {NURBS}.
\newblock {\em International Journal for Numerical Methods in Engineering},
  87:15--47, 2011.

\bibitem{SBVSH11}
M.~A. Scott, M.~J. Borden, C.~V. Verhoosel, T.~W. Sederberg, and T.~J.~R.
  Hughes.
\newblock Isogeometric finite element data structures based on {B}\'ezier
  extraction of {T}-splines.
\newblock {\em International Journal for Numerical Methods in Engineering},
  88:126--156, 2011.

\bibitem{Schillinger2016}
D.~Schillinger, P.~K. Ruthala, and L.~H. Nguyen.
\newblock {L}agrange extraction and projection for {NURBS} basis functions: {A}
  direct link between isogeometric and standard nodal finite element
  formulations.
\newblock {\em International Journal for Numerical Methods in Engineering},
  108(6):515--534, 2016.
\newblock nme.5216.

\bibitem{Kamensky2019}
D.~Kamensky and Y.~Bazilevs.
\newblock {tIGAr}: Automating isogeometric analysis with {FEniCS}.
\newblock {\em Computer Methods in Applied Mechanics and Engineering},
  344:477--498, 2019.

\bibitem{Tirvaudey2020}
M.~Tirvaudey, R.~Bouclier, J.-C. Passieux, and L.~Chamoin.
\newblock Non-invasive implementation of nonlinear isogeometric analysis in an
  industrial {FE} software.
\newblock {\em Engineering Computations}, 37(1):237--261, Jan 2020.

\bibitem{Liu2020}
B.~Liu and D.~Tan.
\newblock A {N}itsche stabilized finite element method for embedded interfaces:
  Application to fluid-structure interaction and rigid-body contact.
\newblock {\em Journal of Computational Physics}, 413:109461, 2020.

\bibitem{Sadeghirad2011}
A.~Sadeghirad, R.~M. Brannon, and J.~Burghardt.
\newblock A convected particle domain interpolation technique to extend
  applicability of the material point method for problems involving massive
  deformations.
\newblock {\em International Journal for Numerical Methods in Engineering},
  86(12):1435--1456, 2011.

\bibitem{Nitsche71}
J.~Nitsche.
\newblock {{\"U}ber ein Variationsprinzip zur L{\"o}sung von
  Dirichlet-Problemen bei Verwendung von Teilr{\"a}umen, die keinen
  Randbedingungen unterworfen sind}.
\newblock {\em Abhandlungen aus dem Mathematischen Seminar der Universit{\"a}t
  Hamburg}, 36:9--15, 1971.

\bibitem{Schillinger2016b}
D.~Schillinger, I.~Harari, M.-C. Hsu, D.~Kamensky, S.~K.~F. Stoter, Y.~Yu, and
  Y.~Zhao.
\newblock The non-symmetric {N}itsche method for the parameter-free imposition
  of weak boundary and coupling conditions in immersed finite elements.
\newblock {\em Computer Methods in Applied Mechanics and Engineering},
  309:625--652, 2016.

\bibitem{Alnaes2014}
M.~S. Aln{\ae}s, A.~Logg, K.~B. {\O}lgaard, M.~E. Rognes, and G.~N. Wells.
\newblock Unified form language: A domain-specific language for weak
  formulations of partial differential equations.
\newblock {\em ACM Trans. Math. Softw.}, 40(2):9:1--9:37, March 2014.

\bibitem{Kirby2006}
R.~C. Kirby and A.~Logg.
\newblock A compiler for variational forms.
\newblock {\em ACM Trans. Math. Softw.}, 32(3):417--444, September 2006.

\bibitem{Logg:2010}
A.~Logg and G.~N. Wells.
\newblock {DOLFIN}: Automated finite element computing.
\newblock {\em ACM Trans. Math. Softw.}, 37(2):20:1--20:28, April 2010.

\bibitem{petsc-web-page}
S.~Balay, S.~Abhyankar, M.~F. Adams, J.~Brown, P.~Brune, K.~Buschelman,
  L.~Dalcin, V.~Eijkhout, W.~D. Gropp, D.~Kaushik, M.~G. Knepley, L.~C.
  McInnes, K.~Rupp, B.~F. Smith, S.~Zampini, and H.~Zhang.
\newblock {PETS}c {W}eb page.
\newblock \url{http://www.mcs.anl.gov/petsc}, 2015.

\bibitem{petsc-user-ref}
S.~Balay, S.~Abhyankar, M.~F. Adams, J.~Brown, P.~Brune, K.~Buschelman,
  L.~Dalcin, V.~Eijkhout, W.~D. Gropp, D.~Kaushik, M.~G. Knepley, L.~C.
  McInnes, K.~Rupp, B.~F. Smith, S.~Zampini, and H.~Zhang.
\newblock {PETS}c users manual.
\newblock Technical Report ANL-95/11 - Revision 3.6, Argonne National
  Laboratory, 2015.

\bibitem{petsc-efficient}
S.~Balay, W.~D. Gropp, L.~C. McInnes, and B.~F. Smith.
\newblock Efficient management of parallelism in object oriented numerical
  software libraries.
\newblock In E.~Arge, A.~M. Bruaset, and H.~P. Langtangen, editors, {\em Modern
  Software Tools in Scientific Computing}, pages 163--202. Birkh{\"{a}}user
  Press, 1997.

\bibitem{Noel2022}
L.~Noel, M.~Schmidt, K.~Doble, J.~A. Evans, and K.~Maute.
\newblock \text{XIGA: An eXtended IsoGeometric Analysis} approach for
  multi-material problems.
\newblock {\em Computational Mechanics}, 2022.

\bibitem{StrFix73}
G.~Strang and G.~J. Fix.
\newblock {\em An Analysis of the Finite Element Method}.
\newblock Prentice-Hall, Englewood Cliffs, NJ, 1973.

\bibitem{Lubarda2020}
M.~V. Lubarda and V.~A. Lubarda.
\newblock {\em Intermediate Solid Mechanics}.
\newblock Cambridge University Press, 2020.

\bibitem{Bazilevs2007a}
Y.~Bazilevs, V.~M. Calo, J.~A. Cottrell, T.~J.~R. Hughes, A.~Reali, and
  G.~Scovazzi.
\newblock Variational multiscale residual-based turbulence modeling for large
  eddy simulation of incompressible flows.
\newblock {\em Computer Methods in Applied Mechanics and Engineering},
  197(1):173--201, 2007.

\bibitem{Bazilevs07c}
Y.~Bazilevs and T.~J.~R. Hughes.
\newblock Weak imposition of {D}irichlet boundary conditions in fluid
  mechanics.
\newblock {\em Computers and Fluids}, 36:12--26, 2007.

\bibitem{Bazilevs08a}
Y.~Bazilevs, V.~M. Calo, T.~J.~R. Hughes, and Y.~Zhang.
\newblock Isogeometric fluid--structure interaction: \text{Theory}, algorithms,
  and computations.
\newblock {\em Computational Mechanics}, 43:3--37, 2008.

\bibitem{Evans2011}
J.~A. Evans.
\newblock {\em Divergence-free {B}-spline Discretizations for Viscous
  Incompressible Flows}.
\newblock {Ph.D.} thesis, University of Texas at Austin, Austin, Texas, United
  States, 2011.

\bibitem{Hsu2010}
M.~C. Hsu, Y.~Bazilevs, V.~M. Calo, T.~E. Tezduyar, and T.~J~R. Hughes.
\newblock Improving stability of stabilized and multiscale formulations in flow
  simulations at small time steps.
\newblock {\em Computer Methods in Applied Mechanics and Engineering},
  199(13):828--840, 2010.
\newblock Turbulence Modeling for Large Eddy Simulations.

\bibitem{Hughes1986}
T.~J.~R. Hughes, L.~P. Franca, and M.~Balestra.
\newblock A new finite element formulation for computational fluid dynamics:
  {V}. {C}ircumventing the {B}abu\v{s}ka--{B}rezzi condition: {A} stable
  {P}etrov--{G}alerkin formulation of the {S}tokes problem accommodating
  equal-order interpolations.
\newblock {\em Computer Methods in Applied Mechanics and Engineering},
  59(1):85--99, 1986.

\bibitem{Franca1992}
L.~P. Franca and S.~L. Frey.
\newblock Stabilized finite element methods: {II}. {T}he incompressible
  {N}avier--{S}tokes equations.
\newblock {\em Computer Methods in Applied Mechanics and Engineering},
  99(2):209--233, 1992.

\bibitem{Kiendl2009}
J.~Kiendl, K.-U. Bletzinger, J.~Linhard, and R.~W\"uchner.
\newblock Isogeometric shell analysis with {K}irchhoff--{L}ove elements.
\newblock {\em Computer Methods in Applied Mechanics and Engineering},
  198:3902--3914, 2009.

\bibitem{Kiendl2011}
J.~Kiendl.
\newblock {\em Isogeometric Analysis and Shape Optimal Design of Shell
  Structures}.
\newblock PhD thesis, Lehrstuhl f\"ur Statik, Technische Universit\"at
  M\"unchen, 2011.

\bibitem{Kiendl2015}
J.~Kiendl, M.-C. Hsu, M.~C.~H. Wu, and A.~Reali.
\newblock Isogeometric {K}irchhoff--{L}ove shell formulations for general
  hyperelastic materials.
\newblock {\em Computer Methods in Applied Mechanics and Engineering},
  291(0):280--303, 2015.

\bibitem{Morganti:2015fy}
S.~Morganti, F.~Auricchio, D.~J. Benson, F.~I. Gambarin, S.~Hartmann, T.~J.~R.
  Hughes, and A.~Reali.
\newblock {Patient-specific isogeometric structural analysis of aortic valve
  closure}.
\newblock {\em Computer Methods in Applied Mechanics and Engineering},
  284:508--520, 2015.

\bibitem{Seder03a}
T.W. Sederberg, J.~Zheng, A.~Bakenov, and A.~Nasri.
\newblock T-splines and {T-NURCCS}.
\newblock {\em ACM Transactions on Graphics}, 22(3):477--484, 2003.

\bibitem{BazHsu12b}
Y.~Bazilevs, M.-C. Hsu, and M.~A. Scott.
\newblock Isogeometric fluid--structure interaction analysis with emphasis on
  non-matching discretizations, and with application to wind turbines.
\newblock {\em Computer Methods in Applied Mechanics and Engineering},
  249--252:28--41, 2012.

\bibitem{Hsu2015}
M.-C. Hsu, D.~Kamensky, F.~Xu, J.~Kiendl, C.~Wang, M.~C.~H. Wu, J.~Mineroff,
  A.~Reali, Y.~Bazilevs, and M.~S. Sacks.
\newblock Dynamic and fluid--structure interaction simulations of bioprosthetic
  heart valves using parametric design with {T}-splines and {F}ung-type
  material models.
\newblock {\em Computational Mechanics}, pages 1--15, 2015.

\bibitem{Casquero2020}
H.~Casquero, X.~Wei, D.~Toshniwal, A.~Li, T.~J.~R. Hughes, J.~Kiendl, and Y.~J.
  Zhang.
\newblock Seamless integration of design and {K}irchhoff--{L}ove shell analysis
  using analysis-suitable unstructured {T}-splines.
\newblock {\em Computer Methods in Applied Mechanics and Engineering},
  360:112765, 2020.

\bibitem{Guo2018}
Y.~Guo, J.~Heller, T.~J.~R. Hughes, M.~Ruess, and D.~Schillinger.
\newblock Variationally consistent isogeometric analysis of trimmed thin shells
  at finite deformations, based on the {STEP} exchange format.
\newblock {\em Computer Methods in Applied Mechanics and Engineering},
  336:39--79, 2018.

\bibitem{Coradello2020}
L.~Coradello, D.~D'Angella, M.~Carraturo, J.~Kiendl, S.~Kollmannsberger,
  E.~Rank, and A.~Reali.
\newblock Hierarchically refined isogeometric analysis of trimmed shells.
\newblock {\em Computational Mechanics}, 66(2):431--447, Aug 2020.

\bibitem{Coradello2021}
L.~Coradello, J.~Kiendl, and A.~Buffa.
\newblock Coupling of non-conforming trimmed isogeometric {K}irchhoff--{L}ove
  shells via a projected super-penalty approach.
\newblock {\em Computer Methods in Applied Mechanics and Engineering},
  387:114187, 2021.

\bibitem{shnapr-code}
\url{https://github.com/david-kamensky/ShNAPr}.
\newblock {ShNAPr} source code.

\bibitem{Kamensky2021}
D.~Kamensky.
\newblock Open-source immersogeometric analysis of fluid--structure interaction
  using {FEniCS} and {tIGAr}.
\newblock {\em Computers \& Mathematics with Applications}, 81:634--648, 2021.
\newblock Development and Application of Open-source Software for Problems with
  Numerical PDEs.

\bibitem{herrema2019penalty}
A.~J. Herrema, E.~L. Johnson, D.~Proserpio, M.~C.~H. Wu, J.~Kiendl, and M.-C.
  Hsu.
\newblock Penalty coupling of non-matching isogeometric {K}irchhoff--{L}ove
  shell patches with application to composite wind turbine blades.
\newblock {\em Computer Methods in Applied Mechanics and Engineering},
  346:810--840, 2019.

\bibitem{Benzaken2021}
J.~Benzaken, J.~A. Evans, S.~F. McCormick, and R.~Tamstorf.
\newblock {N}itsche's method for linear {K}irchhoff--{L}ove shells:
  Formulation, error analysis, and verification.
\newblock {\em Computer Methods in Applied Mechanics and Engineering},
  374:113544, 2021.

\bibitem{Johnson2020}
E.~L. Johnson and M.-C. Hsu.
\newblock Isogeometric analysis of ice accretion on wind turbine blades.
\newblock {\em Computational Mechanics}, 66(2):311--322, Aug 2020.

\bibitem{Liu2021}
N.~Liu, E.~L. Johnson, M.~R. Rajanna, J.~Lua, N.~Phan, and M.-C. Hsu.
\newblock Blended isogeometric {K}irchhoff--{L}ove and continuum shells.
\newblock {\em Computer Methods in Applied Mechanics and Engineering},
  385:114005, 2021.

\bibitem{Zhao2022}
H.~Zhao, X.~Liu, A.~H. Fletcher, R.~Xiang, J.~T. Hwang, and D.~Kamensky.
\newblock An open-source framework for coupling non-matching isogeometric
  shells with application to aerospace structures.
\newblock {\em Computers $\&$ Mathematics with Applications}, 111:109--123,
  2022.

\bibitem{DePrenter2017}
F.~{de Prenter}, C.~V. Verhoosel, G.~J. {van Zwieten}, and E.~H. {van
  Brummelen}.
\newblock Condition number analysis and preconditioning of the finite cell
  method.
\newblock {\em Computer Methods in Applied Mechanics and Engineering},
  316:297--327, 2017.
\newblock Special Issue on Isogeometric Analysis: Progress and Challenges.

\bibitem{Kummer2016}
F.~Kummer.
\newblock Extended discontinuous {G}alerkin methods for two-phase flows: {T}he
  spatial discretization.
\newblock {\em International Journal for Numerical Methods in Engineering},
  109(2):259--289, 2017.

\bibitem{Burman2010}
E.~Burman.
\newblock Ghost penalty.
\newblock {\em Comptes Rendus Mathematique}, 348(21):1217--1220, 2010.

\bibitem{Jomo2019}
J.~N. Jomo, F.~{de Prenter}, M.~Elhaddad, D.~{D'Angella}, C.~V. Verhoosel,
  S.~Kollmannsberger, J.~S. Kirschke, V.~N\"{u}bel, E.~H. {van Brummelen}, and
  E.~Rank.
\newblock Robust and parallel scalable iterative solutions for large-scale
  finite cell analyses.
\newblock {\em Finite Elements in Analysis and Design}, 163:14--30, 2019.

\bibitem{interp-fea-code}
\url{https://github.com/jefromm/interpolation-based-immersed-fea}.
\newblock Source code for {FEniCS}-based solver.

\end{thebibliography}
\end{document}